# Регуляризация линейных задач машинного обучения

*S. Liu[1], S.I. Kabanikhin[2], S.V. Strijhak [3]*
*[1] Novosibirsk State University, Novosibirsk, Russian Federation*
*ORCID: email: liushusng1118@gmail.com*
*[2] Sobolev Institute of Mathematics*
*Siberian Branch of the Russian Academy of Sciences, Novosibirsk, Russian Federation*
*ORCID: email: kabanikh@math.nsc.ru*
*[3] V.P. Ivannikov Institute of System Programming of the Russian Academy of Sciences, Moscow, Russian FederationORCID: email: strijhak@yandex.ru*

**ОГЛАВЛЕНИЕ**

**Аннотация**

В этой работе рассмотрен простейший вариант линейной нейронной сети (ЛНС). Предполагая, что для обучения (построения оптимальной матрицы весов $Q$) мы располагаем набором обучающих пар, то есть нам известны входные данные

$$G = \{g^{(1)}, g^{(2)}, \cdots, g^{(K)}\},$$

а также правильные ответы на эти входные данные

$$H = \{h^{(1)}, h^{(2)}, \cdots, h^{(K)}\}.$$

Исследуем возможности построения матрицы весов $Q$ нейронной сети, которая будет давать правильные ответы на произвольные входные данные, основываясь на связи указанной задачи с системой линейных алгебраических уравнений (СЛАУ). Рассмотрим класс нейронных сетей, в которых каждый нейрон имеет только один выходной сигнал и выполняет линейные операции. Покажем, как подобные ЛНС сводится к СЛАУ.

В силу того, что вопросы $G$ и правильные ответы $H$ нам известны, искомая матрица весов $Q$ должна удовлетворять уравнениям

$$Qg^{(k)} = h^{(k)}, k=1,2,\cdots,K.$$

Требуется восстановить $Q$. В общем случае матрица $Q$ прямоугольная $Q = Q_{MN} = \{q_{mn}\}$, $m$ – номер строки, и $g^{(k)} \in \mathbb{R}^N$, $h^{(k)} \in \mathbb{R}^M$.

Пусть $G_{NK}$ - матрица, составленная из столбцов $g^{(1)}, g^{(2)}, \cdots, g^{(k)}$, а $H_{MK}$– из столбцов $h^{(1)}, h^{(2)}, \cdots, h^{(k)}$. Тогда относительно $Q_{MN}$ получим матричную СЛАУ. $Q_{MN}G_{NK} = H_{MK}$. В этой работе будут изложены методы регуляризации построенной системы.

## 1. Введение

Машинное обучение (МО) - класс методов искусственного интеллекта, которые вместо непосредственного решения задач учатся, применяя решения ряда схожих задач. Для построения такого рода методов используются различные математические и статистические инструменты, методы оптимизации и другие (Сулейманов, 2019) [1].

Мы рассмотрим простейший вариант МО, предполагая, что для обучения нейронных сетей (построения оптимальной матрицы весов $Q$) мы располагаем набором обучающих пар, то есть нам известны входные данные

$$G = \{g^{(1)}, g^{(2)}, \cdots, g^{(K)}\},$$

а также правильные ответы нейронных сетей

$$H = \{h^{(1)}, h^{(2)}, \cdots, h^{(K)}\}.$$

На основе этого набора мы рассмотрим задачу оптимального построения матрицы весов $Q$ нейронной сети, которая будет теперь уже на произвольные входные данные давать правильный ответ.

В данной работе мы рассмотрим некоторые взаимосвязи между исследованием ЛНС и классической теорией и методами регуляризации некорректных задач.

В **разделе 2** ЛНС сводится к системе линейных алгебраических уравнений.

В **разделе 3** строятся примеры ЛНС.

В **разделе 4** приводится краткий обзор методов регуляризации СЛАУ, включая градиентные методы и регуляризацию сдвигом

$$\alpha q_N + A_{KN} q_N = f_N.$$

В **разделе 5** намечены пути обоснования некоторых методов регуляризации ЛНС, включая вопросы согласования ошибки задания входных данных с параметрами регуляризации (снижение эффекта переобучения) и использование априорной и информации (уменьшение размеров системы и ускорение сходимости).

## 2. Сведение ЛНС и СЛАУ

Линейная нейронная сеть — это класс нейронных сетей, в которых каждый нейрон имеет только один выходной сигнал и выполняет линейные операции. В линейной нейронной сети каждый нейрон получает входные данные, умножает их с соответствующими весами, а затем отправляет результат на выход нейрона, суммируя полученные данные (Fukumizu, 1997) [2].

Покажем, как линейная задача машинного обучения сводится к СЛАУ.

Пусть задан набор обучающих пар $\{G, H\}$, то есть, Даны два множества векторов (Kabanikhin, 2023) [3].

Входные $G = \left\{g^{(1)}, g^{(2)}, \cdots, g^{(K)}\right\}$, и выходные $H = \left\{h^{(1)}, h^{(2)}, \cdots, h^{(K)}\right\}$, и также матрица весов $Q$,

$$Qg^{(k)} = h^{(k)}, k = 1,2, \cdots, K. \tag{2.1}$$

Требуется восстановить $Q$.

Матрица $Q$ прямоугольная $Q = Q_{MN} = \{q_{mn}\}$, $m$ – номер строки.

В этом случае $g^{(k)} \in \mathbb{R}^N$, $h^{(k)} \in \mathbb{R}^M$.

Пусть $G_{NK}$ - матрица, составленная из столбцов $g^{(1)}, g^{(2)}, \cdots, g^{(K)}$, а $H_{MK}$ – из столбцов $h^{(1)}, h^{(2)}, \cdots, h^{(K)}$.

Здесь

$$G_{NK} = \begin{pmatrix} g_1^{(1)} & g_1^{(2)} & \cdots & g_1^{(K)} \\ g_2^{(1)} & g_2^{(2)} & & g_2^{(K)} \\ \vdots & & \ddots & \vdots \\ g_N^{(1)} & g_N^{(2)} & \cdots & g_N^{(K)} \end{pmatrix}, H_{MK} = \begin{pmatrix} h_1^{(1)} & h_1^{(2)} & \cdots & h_1^{(K)} \\ h_2^{(1)} & h_2^{(2)} & & h_2^{(K)} \\ \vdots & & \ddots & \vdots \\ h_M^{(1)} & g_M^{(2)} & \cdots & g_M^{(K)} \end{pmatrix}.$$

Тогда относительно $Q_{MN}$ получим матричную систему уравнений.

Здесь

$$Q_{MN} = \begin{pmatrix} q_{11} & q_{12} & \cdots & q_{1N} \\ q_{21} & q_{22} & & q_{2N} \\ \vdots & & \ddots & \vdots \\ q_{M1} & q_{M1} & \cdots & q_{MN} \end{pmatrix}.$$

$$Q_{MN}G_{NK} = H_{MK}. \tag{2.2}$$

Если $M = N = K$, а матрица $G_{NK}$ хорошо обусловлена, то система решается просто

$$Q_{KK} = H_{KK}{G_{KK}}^{-1}. \tag{2.3}$$

**Лемма 2.1**. Если $M = N = K$, определитель матрицы $G_{kK}$ не равен нулю, тогда матрица весов $Q$ определяется однозначно.

Пусть теперь $\det G_{KK} = 0$.

По теореме о сингулярном разложении $G_{KK} = U_{KK}\Sigma_{KK}V_{KK}^*$

Здесь $U_{KK}^*$ и $V_{KK}$ ортогональны, а $\Sigma_{KK} = \mathrm{diag}\left\{\sigma_1, \sigma_2, \cdots, \underbrace{\sigma_{min}}_{\rho}, 0, \cdots, 0\right\}$.

**Лемма 2.2.** Из системы $Q_{KK}G_{KK} = H_{KK}$ мы можем определить $\rho$ линейных комбинаций элементов матриц весов $Q$.

**Доказательство.** Подставим в систему $Q_{KK}G_{KK} = H_{KK}$, сингулярное представление матрицы $G_{KK} = U_{KK}\Sigma_{KK}V_{KK}^*$,

Получим $Q_{KK}U_{KK}\Sigma_{KK}\,V_{KK}^{*} = H_{KK}$

Умножим справа на матрицу $V_{KK}$ получим

$Q_{KK}U_{KK}\Sigma_{KK} \ = H_{KK}V_{KK}$,

$q_{11}u_{11} + q_{21}u_{12}+,\cdots,+q_{K1}u_{1K} = \sigma_1(H_{KK}V_{KK})_{11}$,

$q_{12}u_{21} + q_{22}u_{22}+,\cdots,+q_{K2}u_{2K} = \sigma_2(H_{KK}V_{KK})_{22}$,

$$\cdots\cdots$$

$q_{1k}u_{k1} + q_{2k}u_{k2}+,\cdots,+q_{kk}u_{kk} = \sigma_\rho(H_{kk}V_{kk})_{kk}$,

$$\cdots\cdots$$

$q_{1K}u_{K1} + q_{2K}u_{K2}+,\cdots,+q_{KK}u_{KK} = 0.$

$\rho < K$, является не доопределённой, но мы можем найти нормальное решение.

Для изучение оставшихся случаев отметим, что из (2.2) несложно получить следствие.

$$G_{NK}{}^{\mathrm{T}}Q_{MN}{}^{T} = \ H_{MK}{}^{T}. \tag{2.4}$$

Обозначая $G_{NK}{}^{\mathrm{T}} = A_{KN}$, $H_{MK}{}^{\mathrm{T}} = F_{KM}$, перепишем (2.4) в виде

$$A_{KN}Q_{MN}{}^{T} = F_{KM}. \tag{2.5}$$

Исследуем сначала отдельно СЛАУ для каждого столбца $q^{(m)} = (q_{m1}, q_{m2}, \cdots, q_{mN})^{\mathrm{T}}$, $(m = 1,2,\cdots,M)$, матрицы $Q_{MN}{}^{T}$.

$$A_{KN}q = f. \tag{2.6}$$

Очевидно, что что вопросы разрешимости всех этих СЛАУ, зависят в основном от матрицы $A_{KN}$, поэтому пока что опустим лишние обозначении и рассмотрим систему

$$A_{KN}q = f\ . \tag{2.7}$$

По теореме о сингулярном разложении $A_{KK} = U_{KK}\Sigma_{KK}\,V_{KK}^{*}$

Здесь $U_{KK}^{*}$ и $V_{KK}$ ортогональны, а $\Sigma_{KK} = \mathrm{diag}\left\{\sigma_1, \sigma_2, \cdots, \underbrace{\sigma_{min\{K,N\}}}_{\rho}, 0, \cdots, 0\right\}$.

Следовательно

$$\Sigma_{KK}\,V_{KK}^{*}q = U_{KK}^{*}f \tag{2.8}$$

Величина $\sigma_\rho$ определяет степень устойчивости определения (комбинаций) параметров (матрицы $Q$) вектора $q$. Если $\sigma_\rho > 0$, то

$$(V_{KK}^{*}q_K)_j = \frac{1}{\sigma_j}(U_{KK}^{*}f_K)_j,\ j = 1,2,\cdots,\rho. \tag{2.9}$$

**Лемма 2.3.** *Пусть $1 \le k_0 \le \rho$ – наибольшее целое число при $j = \overline{1,k_0}$, тогда относительно решение LNN можно определить $k_0$ комбинаций (2.8) матрицы $Q$.*

При изучении SLAE возникают следующие вопросы:

1. При каких условиях на $K,\ M,\ N$ и $\{G, H\}$,

LNN имеет единственное устойчивое решение? Другими словами, как формулируются условия корректности LNN?

2. Если условия корректность нарушены, какие методы регуляризации можно применить?

3. Если числа $K$, $M$, $N$ достаточно велики, то решение даже корректной LNN порождает необходимость работы с большими данными, что опять же требует определенной регуляризации (стохастической тензорной и т. п.).

4. В какой мере могут быть объяснены (а в некоторых случаях обоснованы) алгоритмы регуляризации LNN, основанные на применении нейронных сетях?

В данной работе мы постараемся предложить ответы некоторые из вопросов 1–4. Точнее, некоторые ответы на вопросы 1–4.

Разумеется, (2.3) может иметь смысл только в случае, если $G^{-1}$ существует. В противном случае задача (2.2) некорректна и появляется необходимость строить регуляризирующие алгоритмы.

**3. Примеры ЛНС**

Простейший пример линейные нейронные сети.

Пусть

$h = \sum_{i=1}^{N} g_i q_i + b\ ,\ i = 1,2,\cdots N.\ y = f(h)\ .$ (3.1)

где $h$ – выходной вектор, $q$ – вектор весов, $g$ – входной вектор, $b$ – параметр смещения, $f$ – функция активации (Уоссермен,1992) [4].

Нейрон представляет собой своеобразную единицу обработки информации в нейронной сети.

1. Набор синапсов, которые характеризуются своими весами, на которые умножаются входные сигналы.

2. Сумматор, который суммирует входные сигналы, взвешенные относительно соответствующих синапсов нейронов.

3. Функция активации, обеспечивающая сжатое отображение взвешенного суммарного сигнала таким образом, что выходной сигнал имеет ограниченный диапазон $[0,1]$ или$[-1,1]$ (Сысоев, т. п., 2014) [5].

Передаточная функция нейрона линейной нейронной сети является линейной функцией, поэтому на выходе может быть любое значение, а на выходе перцептрона может быть только 0 или 1 (Горбань,1998) [6].

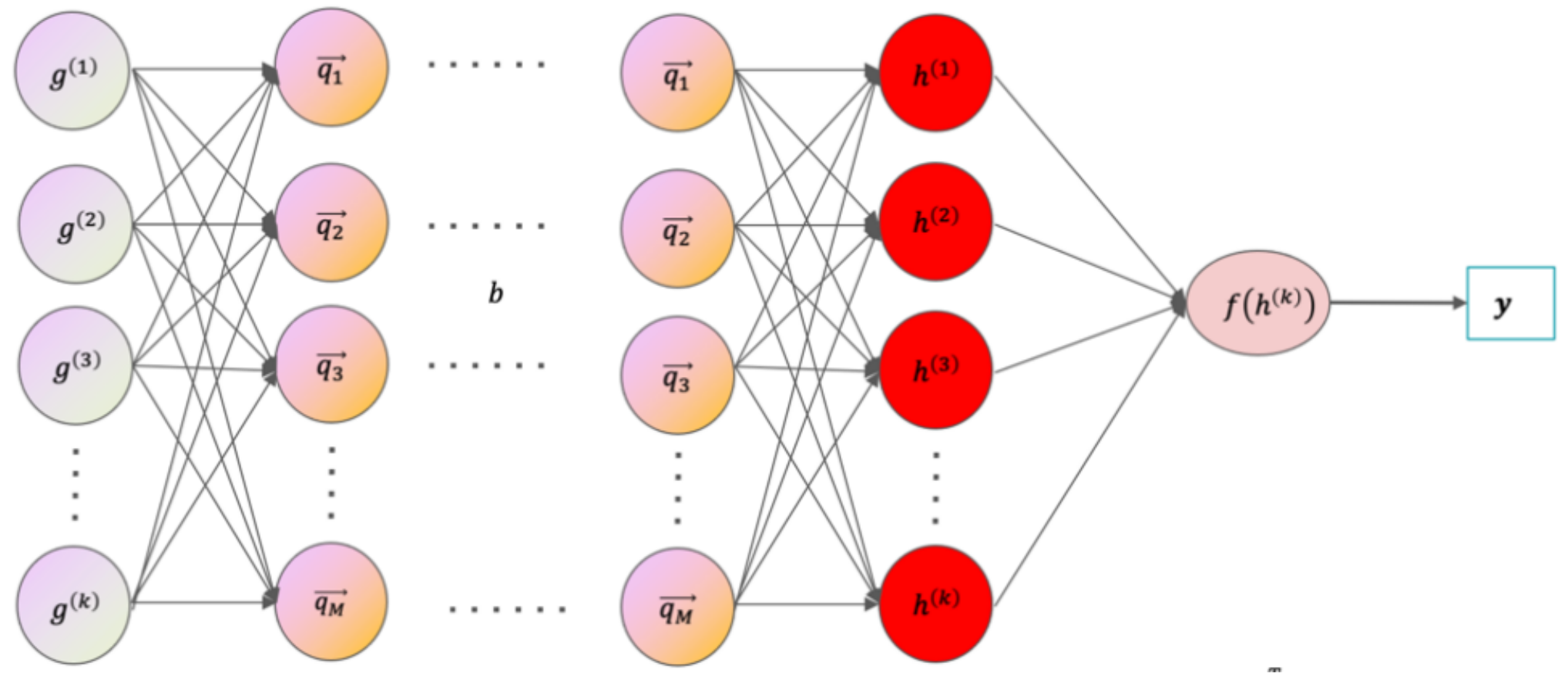

Рис. 1: Структура линейной нейронной сети.

**3.1. Численные эксперименты**

**Пример 3.1.**

Пусть $g^{(k)}$ – результаты анализов пациента, $h^{(k)}$- характеристика состояние здоровья пациента, условно говоря “предварительный диагноз”. Компоненты $g_n^{(k)}$ вектора $g^{(k)}$ принимает значение 1 – если результат анализа с номером $n$ содержатся в допустимых пределах, и значение 0 – если результат анализа выходит за пределы нормы, $n = 1,2,\cdots,N$, $k = 1,2,\cdots,K.$

Компоненты вектора $h^{(k)}$ принимают значение:

1– если пациент здоров;

2. – если он “скоро здоров, чем болен”;

3. – если он “скоро болен, чем здоров”;

4. – если пациент болен.

Таким образом, в нашем примере $g^{(k)} \in \mathbb{R}^2$, $h^{(k)} \in \mathbb{R}^4$.

$G = \begin{bmatrix} 1 & 0 \\ 0 & 0 \end{bmatrix}, H = \begin{bmatrix} 1 \\ 1 \end{bmatrix}$

Ранг Матрицы G: 1

Матрица G не является полно ранговой.

Pseudo Solution:

[[1. 0.]

[0. 0.]]

Final Pseudo Solution:

[[9.8998642e-01 1.8486163e-10]

[9.8998642e-01 1.8486163e-10]]

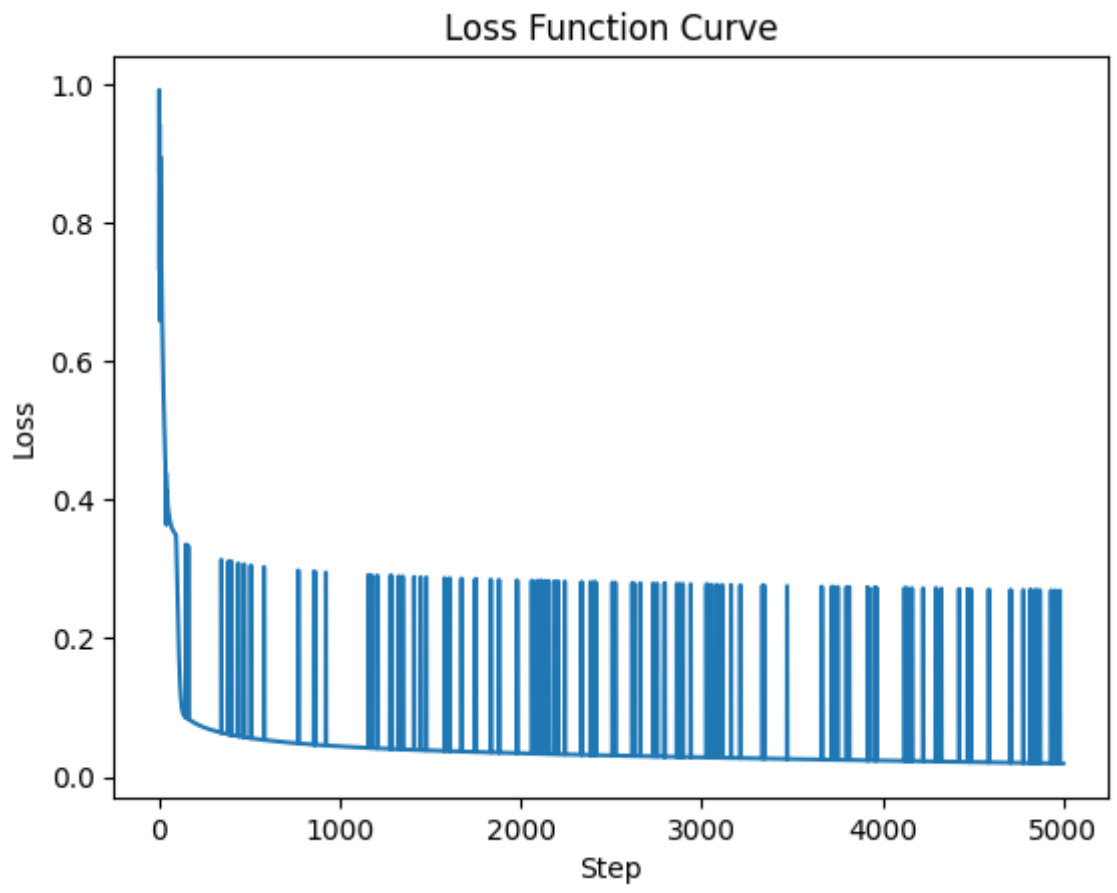

Рис. 2: Вывод значения функции потерь на каждой итерации обучения.

**Пример 3.2.**

Рассмотрим более детальное описание результатов анализов и классификацию оценки состояние здоровья пациентов.

Допустим, что результаты анализов пациента с номером $k = 1,2,\cdots,K$. Принимают значение от 1 (самый благоприятный результат, “норма”) до 20 (самый нежелательный результат), то есть, $g^{(k)} \in \mathbb{R}^{20}$.

Состояние здоровья пациента также упорядочим от 1 (абсолютно здоров) до 30 (серьезно болен), то есть, $h^{(k)} \in \mathbb{R}^{30}$.

В этом случае матрица исконных чёсов $Q_{30,20}$ удовлетворяет матричной метем управлений $Q_{30,20}G_{20,K}$=$H_{30,K}$.

Rank of Matrix $G$: 20

Матрица $G$ является полно ранговой.

Pseudo Solution: 20x20

Final Pseudo Solution: 30x20

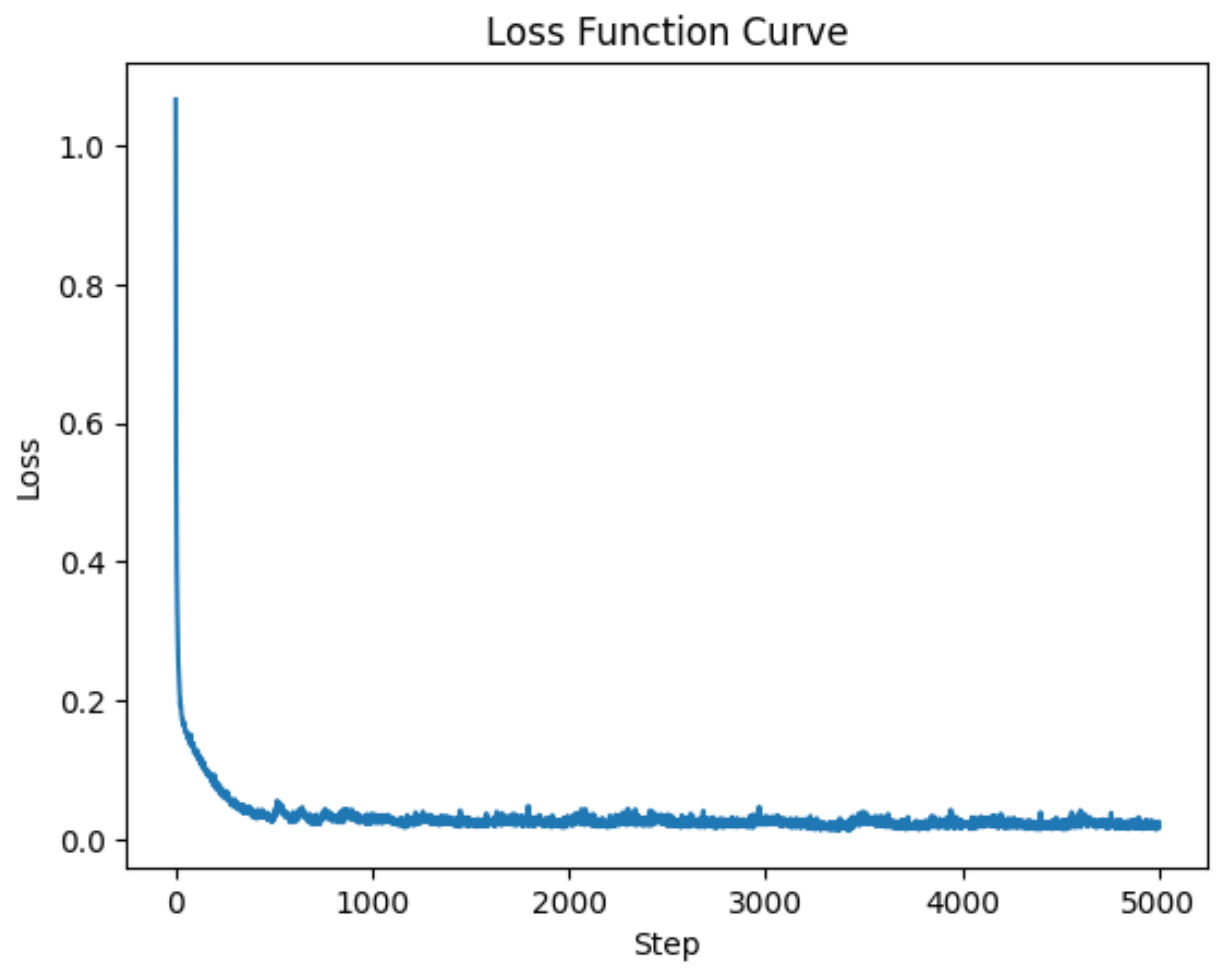

Рис. 3: Вывод значения функции потерь на каждой итерации обучения.

**Пример 3.2.**

Построение модели анализа данных анализа крови на основе линейной нейронной сети:

Импорт библиотек

Сначала импортируем необходимые библиотеки Python: numpy, tensorflow, matplotlib, scipy и sklearn.

Подготовка данных

Соберите и очистите данные анализа крови. Пусть результаты анализов пациентов пронумерованы $k = 1,2,\cdots,K$. В таблице ниже представлены название теста, результат анализа и эталонный стандарт.

Таблица 1. Анализ данных анализа крови.

| Number | Code | Item Name | Result | Reference Range | Unit |
|---|---|---|---|---|---|
| 1 | age | age | | 0-18 | Years |
| 2 | WBC | white blood cell count | | 4-10 | 10*9/L |
| 3 | #LYM | lymphocyte Count | | 0.8-4 | 10*9/L |
| 4 | #MON | Intermediate Cell Count | | 1.1-1.2 | 10*9/L |
| 5 | #GRAN | Central Cell Count | | 2-7 | 10*9/L |
| 6 | % LYM | Lymphocyte Percentage | | 20-40 | % |
| 7 | %MON | Intermediate Cell | | 3-14 | % |
| | | Percentage | | | |
| 8 | %GRAN | Neutrophil Percentage | | 50-70 | % |
| 9 | RBC | Red Blood Cell Count | | 3.5-5.5 | 10*12/L |
| 10 | HGB | Hemoglobin | | 110-160 | g/L |
| 11 | HCT | Hematocrit | | 37-54 | % |
| 12 | MCV | Mean Corpuscular Volume | | 80-100 | fL |
| 13 | MCH | Mean Corpuscular Hemoglobin | | 27-34 | pg |
| 14 | MCHC | Mean Corpuscular Hemoglobin Concentratio | | 320-360 | g/L |
| 15 | RDW-CV | Red Cell Distribution Width - Coefficient of Variation | | 11-16 | % |
| 16 | RDW-SV | Red Cell Distribution Width - Standard Deviation | | 35-56 | fL |
| 17 | PLT | Platelet Count | | 100-300 | 10*9/L |
| 18 | MPV | Mean Platelet Volume | | 6.5-12 | fL |
| 19 | PDW | Platelet Distribution Width | | 9-17 | fL |
| 20 | PCT | Plateletcrit | | 0.108-0.282 | % |

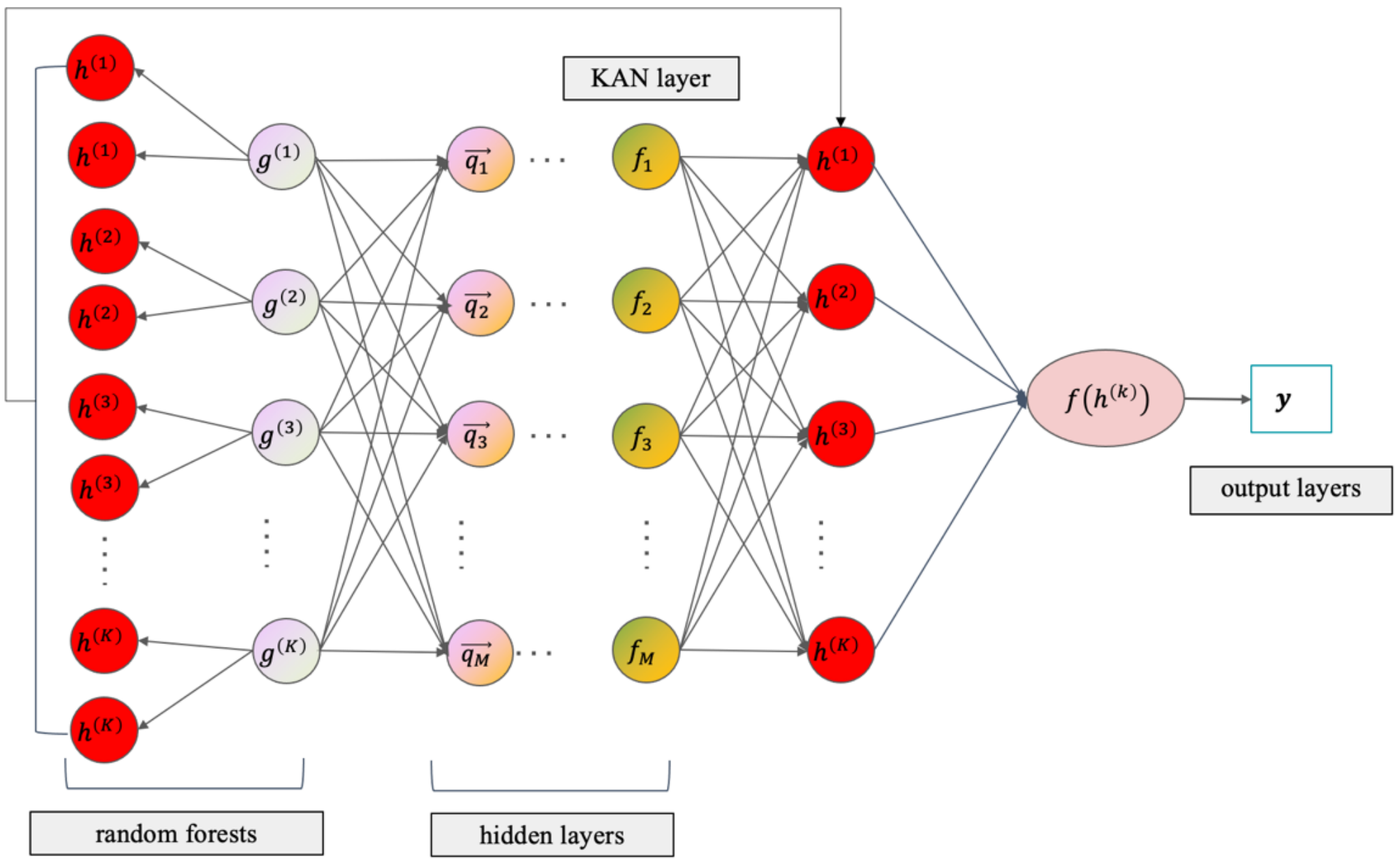

Рис. 4: Структура модели анализа данных анализа крови на основе линейной нейронной сети.

Инженерия признаков

Преобразовать данные анализа крови в признаки, подходящие для линейной нейронной сети. Используйте методы отбора признаков, уменьшения размерности и стандартизации для обработки данных. Пусть $g^{(k)}$ — это результаты обследования пациента, а $h^{(k)}$ — это характеристика здоровья пациента, представляющая собой "предварительный диагноз". Компоненты вектора $g^{(k)}$ обозначаются как $g_n^{(k)}$, где 0 - если результат анализа с номером n находится в пределах допустимых норм, и 1 - если результат анализа выходит за пределы нормы. Этот процесс использует метод случайного леса для классификации результатов обследования пациентов и определения значений компонентов вектора $h^{(k)}$. Масштабный коэффициент данных устанавливается равным 10,9.

Настройка параметра слоев

Указать количество нейронов в каждом слое: layers = [10, 20, 100, 100, 10, 10].

Метод инициализации

Использовать метод инициализации HeNormal для инициализации весов.

Структура сети

Каждый слой включает в себя Dense полносвязный слой, BatchNormalization слой, пользовательский PReLU слой и Dropout слой (dropout rate = 0.01) и слой KAN.

Пользовательский слой PReLU

PReLU (Parametric ReLU) — это функция активации. Класс PReLU наследуется от tf.keras.layers.Layer и реализует методы инициализации, построения и вызова.

Пользовательский слой KANLinear

Класс KANLinear является пользовательским полносвязным слоем с несколькими параметрами и методами, включая reset_parameters и curve2coeff.

Класс NeuralNetworkModel

Класс NeuralNetworkModel наследуется от tf.keras.Model и реализует пользовательскую линейную нейронную сеть, включающую несколько полносвязных слоев, слоев BatchNormalization, слоев PReLU и слоев KANLinear.

Входные данные

$g^{(k)}$ и $h^{(k)}$ - входные данные, параметр используется для различения этапов тренировки и вывода.

Метод прямого распространения

Данные последовательно проходят через слои Dense, BatchNormalization, PReLU и Dropout, а также KAN.

Выходные данные

Выход последнего слоя преобразуется в матрицу формы 20x10, то есть для каждого образца имеется 10 выходных узлов.

Обучение модели

Обучение модели глубокого обучения с использованием подготовленного набора данных и п роверка с использованием техники кросс-валидации. Вычислили ошибку между окончательным решением псевдообратной матрицы и истинным значением

Error $= \|H - QG\| =$1.739243.

Оценка модели

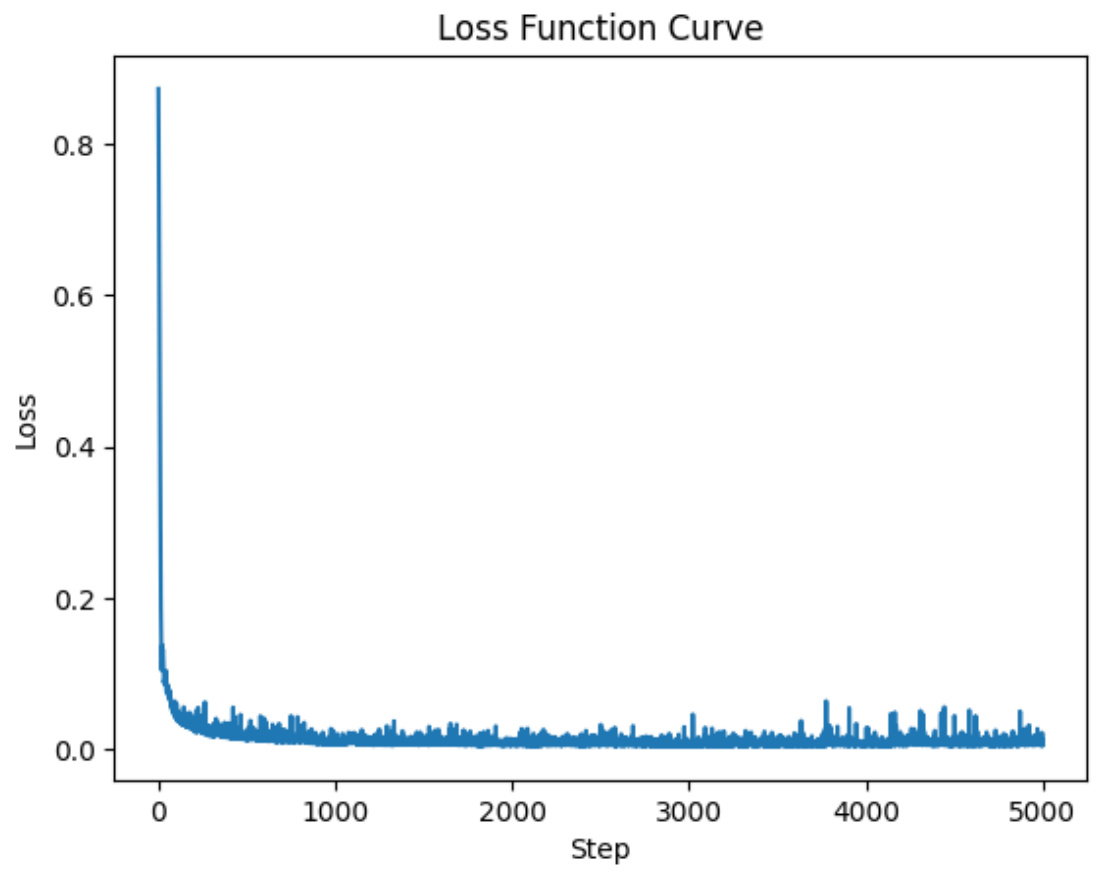

Рис. 5: Вывод значения функции потерь на каждой итерации обучения.

## 4. Обзор регуляризации ЛНС

Обучение линейных нейронных сетей осуществляется путем корректировки весов нейронов, с целью минимизировать ошибку между прогнозируемыми и фактическими значениями в обучающих данных (Kwon, 1992) [7]. Для этой цели используются методы оптимизации, чаще всего градиентный спуск (Duchi, 2011) [8]. Перед обучением необходимо определить функцию потерь, которая измеряет разницу между прогнозируемым значением и фактическим значением. Обычно используется стандартная ошибка или средняя абсолютная погрешность. После определения функции потерь можно приступать к обучению. На каждой итерации градиентного спуска вес нейрона обновляется в направлении, противоположном градиенту функции потерь. Это позволяет свести к минимуму ошибки и повысить предсказательную способность сети. Обучение линейных нейронных сетей может быть довольно быстрым и эффективным, особенно для регрессионных задач с большим количеством функций. Однако для более сложных задач, таких как классификация изображений, требуется более сложная архитектура нейронной сети.

В литературе по нейронным сетям регуляризация понимается как метод, используемый для предотвращения переобучения модели. Эти методы уменьшают влияние шума в данных и делают модель более устойчивой к изменениям входных данных. Мы упомянем только 4 основных метода в этой связи:

1. Регуляризация L1 и L2 (Регуляризация А. Н. Тихонова): штрафное слагаемое добавляется к функции потерь (для L1-регуляризации - абсолютное значение параметра, для L2-регуляризации - квадратичная норма параметра). Этот метод может быть использован для выбора наилучшего значения регуляризованного гиперпараметра. Например, для оценки качества модели при различных значениях гиперпараметра применяется перекрестная проверка, и выбирается значение, дающее наилучшее качество (Kalivas, 2012) [9].

2. Отсев: некоторые нейроны случайным образом отключаются во время обучения, тем самым сокращая переподготовку (Srivastava, 2014) [10].

3. Ранняя остановка: если ошибки в наборе тестовых данных начинают увеличиваться, прекратите обучение (Prechelt, 2002) [11].

4. Пакетная нормализация: нормализуйте данные на каждом уровне нейронной сети, чтобы ускорить обучение и сократить время переобучения (Ioffe and Szegedy) [12].

5. Расширение данных: проблема переобучения нейронных сетей решается путем улучшения обучающих данных, но затраты на улучшение данных очень высоки, и иногда данные даже невозможно улучшить. Затем обучающий набор может быть расширен путем добавления данных (Simard, 2003) [13]. В области компьютерной обработки изображений из-за очень большого размера входного изображения трудно найти так много связанных изображений, что требует использования методов улучшения данных: горизонтальное переворачивание изображений, обрезка изображений, поворот, искаженная нумерация изображений и т. д., которые также могут быть применены к линейным нейронные сети.

Отметим взаимосвязь между этими методами и классическими методами регуляризации.

### 4.1. Методы регуляризации

#### 4.1.1. L1 регуляризация

Регуляризация L1 также называется регрессией лассо, Ω – (штрафная функция) представляет собой сумму абсолютных значений весовых параметров в весовой матрице (Kalivas, 2012) [9]：

$$\Omega(Q) = \|Q\|_1 = \sum_i \sum_j |q_{ij}|, i = 1,2,\cdots N, j = 1,2,\cdots N. \quad (4.1)$$

Мы умножаем штраф на $\alpha$ и добавляем к функции потерь $\mathcal{L}(Q)$.

$$\hat{\mathcal{L}}(Q) = \alpha\|Q\|_1 + \mathcal{L}(Q). \quad (4.2)$$

где $\hat{\mathcal{L}}(Q)$ — это обновленная функция потерь, $\alpha$ — это параметр регуляризации.

### 4.1.2. L2-регуляризация

Регуляризация L2 является наиболее распространенным типом из всех методов регуляризации и также обычно упоминается как уменьшение веса. При регуляризации L2 функция потерь нейронной сети расширяется на слагаемое регуляризации Ω.

$$\Omega(Q) = \|Q\|_2^2 = \sum_i \sum_j q_{ij}^2, i = 1,2,\cdots N, j = 1,2,\cdots N. \quad (4.3)$$

$$\hat{\mathcal{L}}(Q) = \frac{\alpha}{2}\|Q\|_2^2 + \mathcal{L}(Q) = \frac{\alpha}{2}\sum_i \sum_j q_{ij}^2 + \mathcal{L}(Q). \quad (4.4)$$

Члены регуляризации L2 и L1 добавляются к функции потерь $\mathcal{L}(Q)$ из LNN. Роль терминов регуляризации заключается в том, чтобы избежать переобучения. Сначала определяется параметр регуляризации a, и пропорция члена регуляризации в общих потерях регулируется путем корректировки параметра $\alpha$. Затем значения членов регуляризации L2 и L1 были вычислены отдельно. Термин регуляризации L2 вычисляется путем суммирования квадратов всех весов и умножения на 0.5, в то время как термин регуляризации L1 вычисляется путем суммирования абсолютных значений всех весов $Q$. Наконец, условия регуляризации L2 и L1 умножаются на параметры регуляризации $\alpha$ и добавляются к потере $\mathcal{L}(Q)$, так что условия потери и регуляризации будут сведены к минимуму одновременно в процессе обучения, чтобы достичь цели сбалансированной подгонки и улучшения способности к обобщению.

$$\hat{\mathcal{L}}(Q) = \mathcal{L}(Q) + \alpha\|Q\|_1 + \frac{\alpha}{2}\|Q\|_2^2 = \mathcal{L}(Q) + \alpha\sum_i \sum_j |q_{ij}| + \frac{\alpha}{2}\sum_i \sum_j q_{ij}^2. \quad (4.5)$$

### 4.1.3. Dropout

Dropout часто называют «методом отсева» или «методом исключения».

В процессе обучения некоторые нейроны случайным образом отключаются, что позволяет уменьшить переобучение (Ketkar, 2017) [14].

Dropout — это метод регуляризации, при котором некоторые нейроны игнорируются случайным образом. Это означает, что изменения веса не применяются к данному нейрону при обратном проходе и что их влияние на активацию нижележащих нейронов временно стирается при прямом проходе. Веса нейронов внутри нейронной сети занимают свое место в сети по мере ее обучения.

$$h = f\left(\sum_{i=1}^{N} g_i q_i + b\right), i = 1,2,\cdots N. \quad (4.6)$$

Этап обучения:

$$y = h \odot B, B_i \sim Bernoulli(p). \quad (4.7)$$

Этап обучения:

$$y = (1-p)h \odot B. \quad (4.8)$$

где $p$ — это вероятность сохранения активированного нейронного узла.

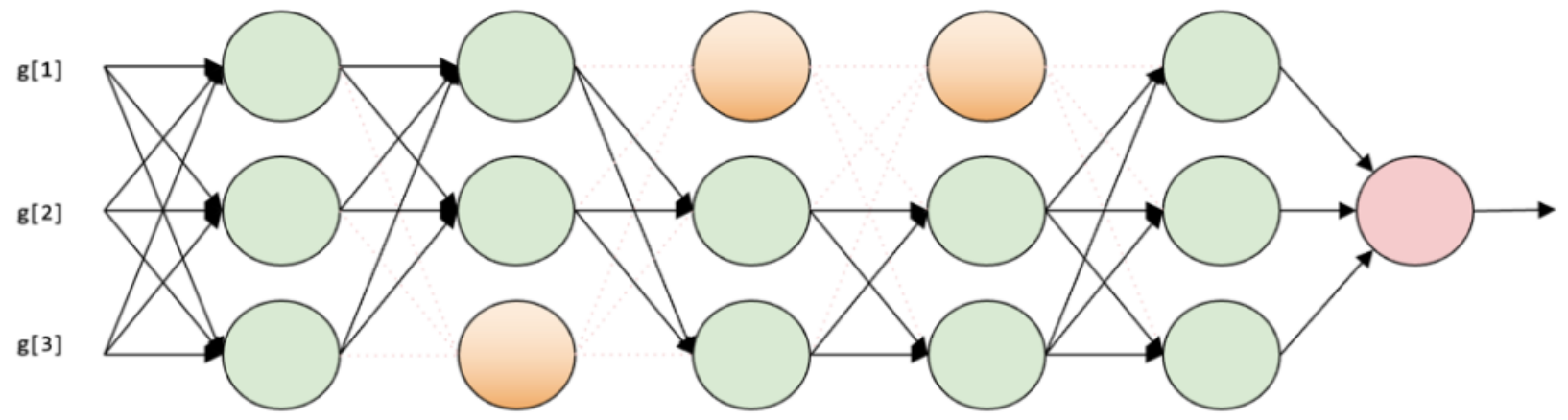

Рис. 6: Структура Dropout.

На рис. 6. представлена непрерывная нейронная сеть с прямой связью, в которой используется Dropout (Wang, 2018) [15].

**4.1.4. Пакетная нормализация**

Пакетная нормализация (Batch normalization)— это другой метод регуляризации может нормализовать набор активаций в слое. Нормализация выполняется путем вычитания пакетного среднего значения из каждой активации (Ioffe and Szegedy, 2015) [12].

Скорость сходимости модели ускоряется, и, что более важно, проблема “градиентной дисперсии” в глубокой сети в определенной степени решается, что упрощает и делает более стабильным обучение модели глубокой сети.

Пакетная нормализация может использоваться в качестве слоя нейронных сетей, размещаемого перед функциями активации. Алгоритм работы пакетной нормализации выглядит так, как показано на рис. 3.

Нормализация данных на каждом слое нейронной сети для ускорения обучения и снижения переобучения (Luo, et. al., 2018) [16].

$$\mu_B = \frac{1}{N}\sum_{i=1}^{N} g_i\,,\ i = 1,2,\cdots N \tag{4.9}$$

$$\sigma_B^2 = \frac{1}{N}\sum_{i=1}^{N}(g_i - \mu_B)^2, \tag{4.10}$$

$$\hat{g}_i = \frac{g_i - \mu_B}{\sqrt{\sigma_B^2 + \epsilon}}., \tag{4.11}$$

$$h_i = \gamma\hat{g}_i + \beta = BN_{\gamma\beta}(g_i). \tag{4.12}$$

Где γ, β — обучаемые параметры.

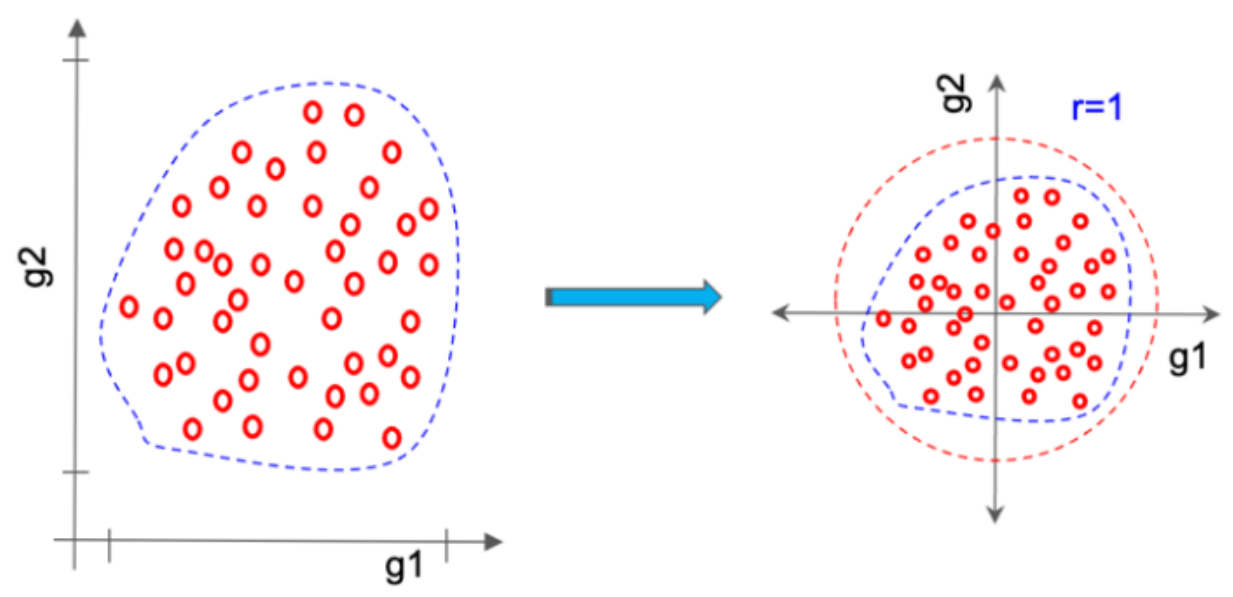

Рис. 7: Структура пакетной нормализации.

**4.1.5. Ранняя остановка**

Ранняя остановка (Early stopping) означает, что мы периодически сохраняем обучаемые параметры и отслеживаем ошибку валидации. После прекращения обучения мы возвращаем обучаемые параметры в точное положение, в котором ошибка валидации начала расти, а не в положение последнего параметра.

Другой способ представить себе раннюю остановку — это очень эффективный алгоритм выбора гиперпараметров, который устанавливает число эпох как абсолютно лучшее. По сути, ограничивает процедуру оптимизации небольшой частью обучаемого пространства параметров, близкой к начальным параметрам.

Можно также показать, что в случае простых линейных моделей с квадратичными функциями ошибок и простым градиентным спуском ранняя остановка эквивалентна $L2$ - регуляризации (Prechelt, 2002) [11].

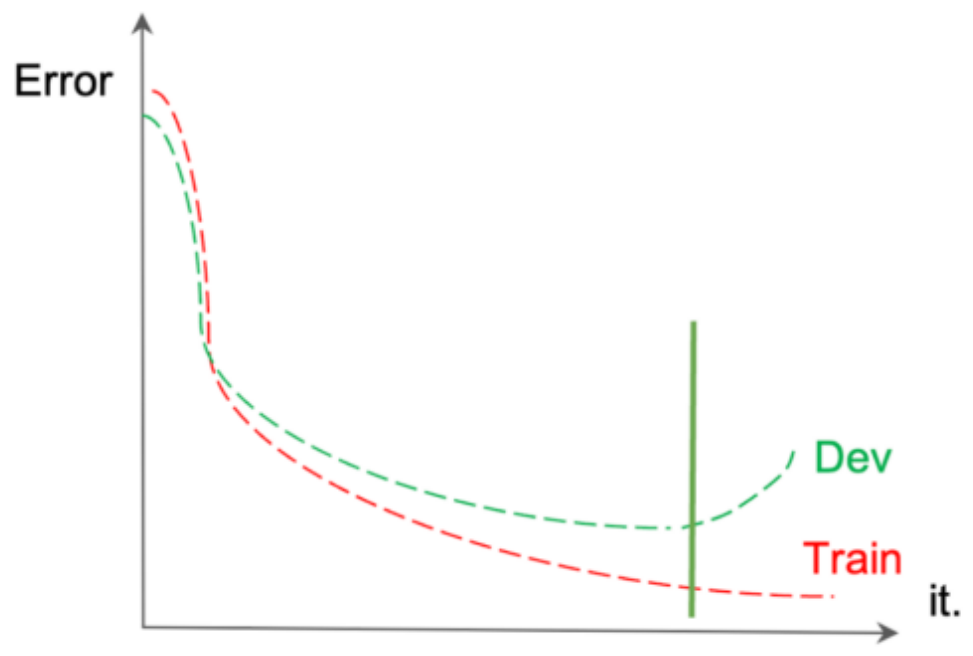

Рис. 8: Ранняя остановка: прекратить обучение, если ошибка на валидационном наборе данных начинает расти.

#### 4.1.6. Расширение данных

Предположим, мы обучаем нейронную сеть для выполнения задачи классификации изображений, надеясь решить проблему переобучения путем усиления обучающих данных, но усиление данных обходится очень дорого, и иногда данные не могут быть усилены. Затем можно увеличить обучающий набор, добавив изображения (Simard, 2003) [13].

Результаты улучшения данных включают в себя предотвращение переобучения; повышение надежности модели и снижение чувствительности модели к изображению; увеличение объема обучающих данных для улучшения обобщающей способности модели; предотвращение неравномерных выборок и уменьшение доли неравномерных выборок за счет применения некоторых методов улучшения данных на небольших выборках.

Распространенными методами улучшения данных в основном являются переворачивание, поворот, обрезка, масштабирование, панорамирование и встряхивание. Наиболее часто используемыми методами преобразования пикселей для обработки изображений являются: добавление резкого шума, гауссовского шума, размытия по Гауссу, регулировка контрастности HSV, регулировка яркости, насыщенности, выравнивание гистограммы, регулировка баланса белого и т. д.

Функция улучшения искажения данных позволяет преобразовывать существующие изображения с сохранением их названий. Это включает в себя такие улучшения, как преобразование геометрии и цвета, случайное стирание, состязательное обучение и перенос нейронного стиля (Shorten, 2019) [17].

#### 4.1.7. Сингулярное разложение

Сингулярное разложение (SVD): позволяет уменьшить размерность данных, сжать информацию, уменьшить количество параметров модели и ускорить обучение. Применение SVD к матрице весов между двумя слоями нейронной сети очень эффективно. Мы можем применить SVD к матрице активации нейронов в конкретном слое и найти наиболее важные сингулярные значения, которые будут соответствовать наиболее важным нейронам. Это позволяет уточнить, какие нейроны и слои вносят наибольший вклад в работу нейронной сети и какие нейроны и слои можно опустить без ущерба для качества модели.

### 4.2. Проблема исчезновения градиента/взрыва

Нейронные сети полагаются на обратное распространение для передачи информации туда и обратно. Однако по мере уменьшения количества слоев градиент часто становится все меньше и меньше, что приводит к почти неизменному весу нижних слоев, и мы не можем прийти к хорошему решению. Это проблема исчезновения градиента. С другой стороны, градиент может становиться все больше и больше, что приводит к отклонению алгоритма, что является проблемой взрыва градиента. Ксавье и Бенджио указали, что популярная в то время функция активации Sigmoid и метод инициализации со средним значением 0 и стандартным отклонением 1 имели серьезные проблемы: дисперсия выходных данных каждого слоя была намного больше дисперсии его входных данных, и увеличение дисперсии было насыщенным в верхний слой (Glorot and Bengio, 2010) [18].

Функция активации может быть разделена на две категории — "функция активации насыщения" и "функция активации ненасыщения".

Sigmoid и Tanh — это "функции активации насыщения". Сигмовидная функция должна сжимать входные данные до диапазона [0,1], а функция Tanh должна сжимать входные данные до диапазона [-1,1].

Метод инициализации функцией активации Sigmoid и средним значением, равным 0, и стандартным отклонением, равным 1, имеет серьезную проблему: дисперсия выходных данных каждого слоя намного больше, чем дисперсия его входных данных. Увеличение дисперсии приводит к насыщению функции активации на верхнем уровне. На самом деле это происходит потому, что среднее значение логической функции равно 0,5 вместо 0 (среднее значение функции Tanh равно 0, которая работает немного лучше). Как показано на рисунке 5, можно видеть, что, когда входные данные функции Sigmoid становятся больше (отрицательными или положительными), насыщенность функции равна 0 или 1, а производная очень близка к 0. Следовательно, когда начинается обратное распространение, у него почти нет градиента для обратного распространения по сети.

Функция ReLU представляет собой "скорректированную линейную единицу измерения”, которая является максимальной функцией входных данных со сверточным изображением. Функция ReLU устанавливает все отрицательные значения в матрице x равными нулю, а остальные значения остаются неизменными. Вычисление функции ReLU выполняется после свертки, поэтому она такая же, как функция Tanh и сигмовидная функция, и относится к "нелинейной функции активации". Этот контент был впервые предложен Джеффом Хинтоном (LeCun, 2015) [19]. Параметризованный линейный блок коррекции (PReLU) можно рассматривать как вариант негерметичного ReLU. В PReLU наклон отрицательной части определяется на основе данных, а не задается заранее. Leaky ReLU — это первый вариант ReLU. Leaky ReLU успешно решает проблему ослабления градиента ReLU, присваивая ненулевой наклон отрицательной части (Maas, 2013) [20]. Впоследствии был предложен параметр ReLU (PReLU), который требует изучения отрицательной части вместо использования заранее определенных значений (He, 2015) [21].

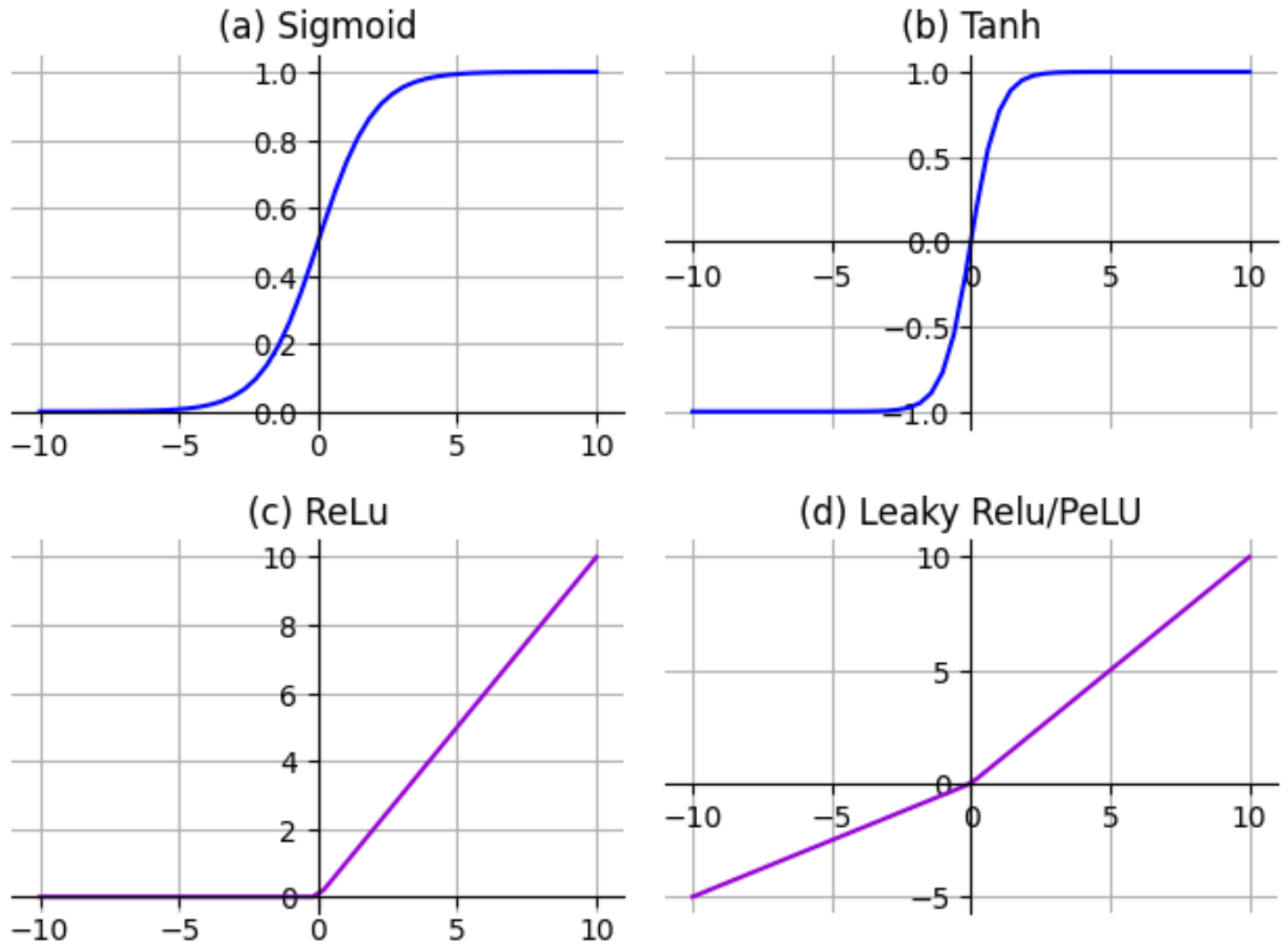

Рис. 9: Функции активации.

Sigmoid: $S(x) = \frac{1}{1+e^{-x}}$ . (4.13)

Tanh: $T(x) = \frac{e^{x}-e^{-x}}{e^{x}+e^{-x}}$. (4.14)

ReLU:$R(x) = \max(0, x)$. (4.15)

Leaky ReLU: $LR(x) = \begin{cases} x_i, & x_i \geq 0 \\ \frac{x_i}{\alpha_i}, & x_i < 0 \end{cases}$. (4.16)

Где $\alpha_i$ — это фиксированный параметр в интервале $(1, +\infty)$. Параметр $\alpha_i$ в PReLU изменяется в зависимости от данных; параметр $\alpha_i$ в Leaky ReLU остается неизменным.

### 4.3. Переоснащение и недооснащение

Модель изучает взаимосвязь между входными данными (называемыми объектами) и выходными данными (называемыми метками) из набора обучающих данных. В процессе обучения модель получает как объекты, так и метки и учится сопоставлять первые со вторыми. Например, предположим, что мы хотим построить линейную нейронную сеть для прогнозирования цен на недвижимость. Введите начальную дату создания данных, площадь участка (квадратные футы) и вместимость гаража, и результатом будет цена продажи дома. Обученная модель оценивается на тестовом наборе. На тестовом наборе мы только предоставляем ей характеристики, и она делает прогнозы. Мы сравниваем прогнозируемые результаты с известными метками тестового набора, чтобы рассчитать точность. Важной частью нашей генерации данных является добавление случайного шума к метке. В любом реальном процессе, будь то естественный или искусственный, данные не полностью подходят для определения тенденции. В отношениях всегда присутствует шум или другие переменные, которые мы не можем измерить. В случае цен на жилье, из-за других факторов, влияющих на цены на жилье, цены не полностью совпадают.

• Переобучение: чрезмерная зависимость от обучающих данных

• Недостаточная подгонка: невозможно изучить взаимосвязь в обучающих данных

• Высокая дисперсия: модель существенно изменяется в зависимости от обучающих данных.

• Высокое отклонение: предположения о модели приводят к игнорированию обучающих данных.

• Чрезмерная и недостаточная подгонка приводят к ухудшению способности к обобщению в тестовом наборе.

• Проверочный набор, используемый для настройки модели, может предотвратить недостаточную и чрезмерную подгонку (Van der Aalst, 2010) [22], (Koehrsen, 2018) [23].

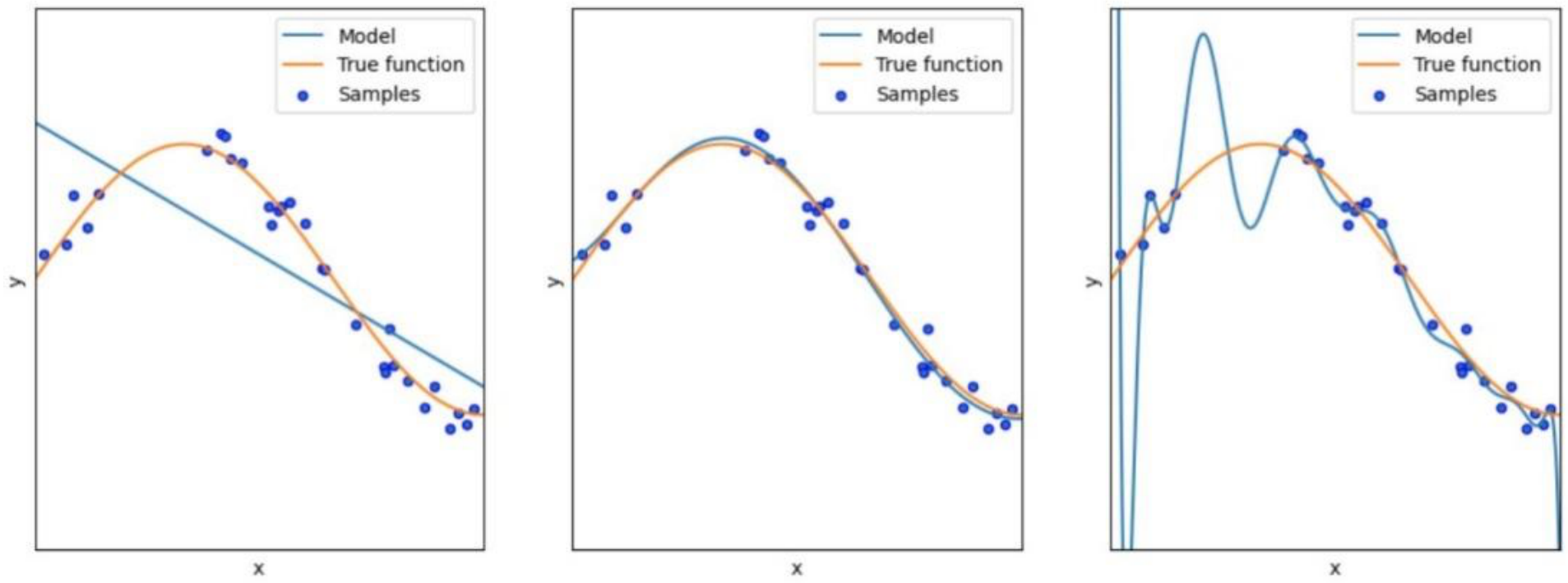

Рис. 10: Train-Dev-Test в сравнении с подгонкой модели.

## 5. Обзор регуляризации СЛАУ

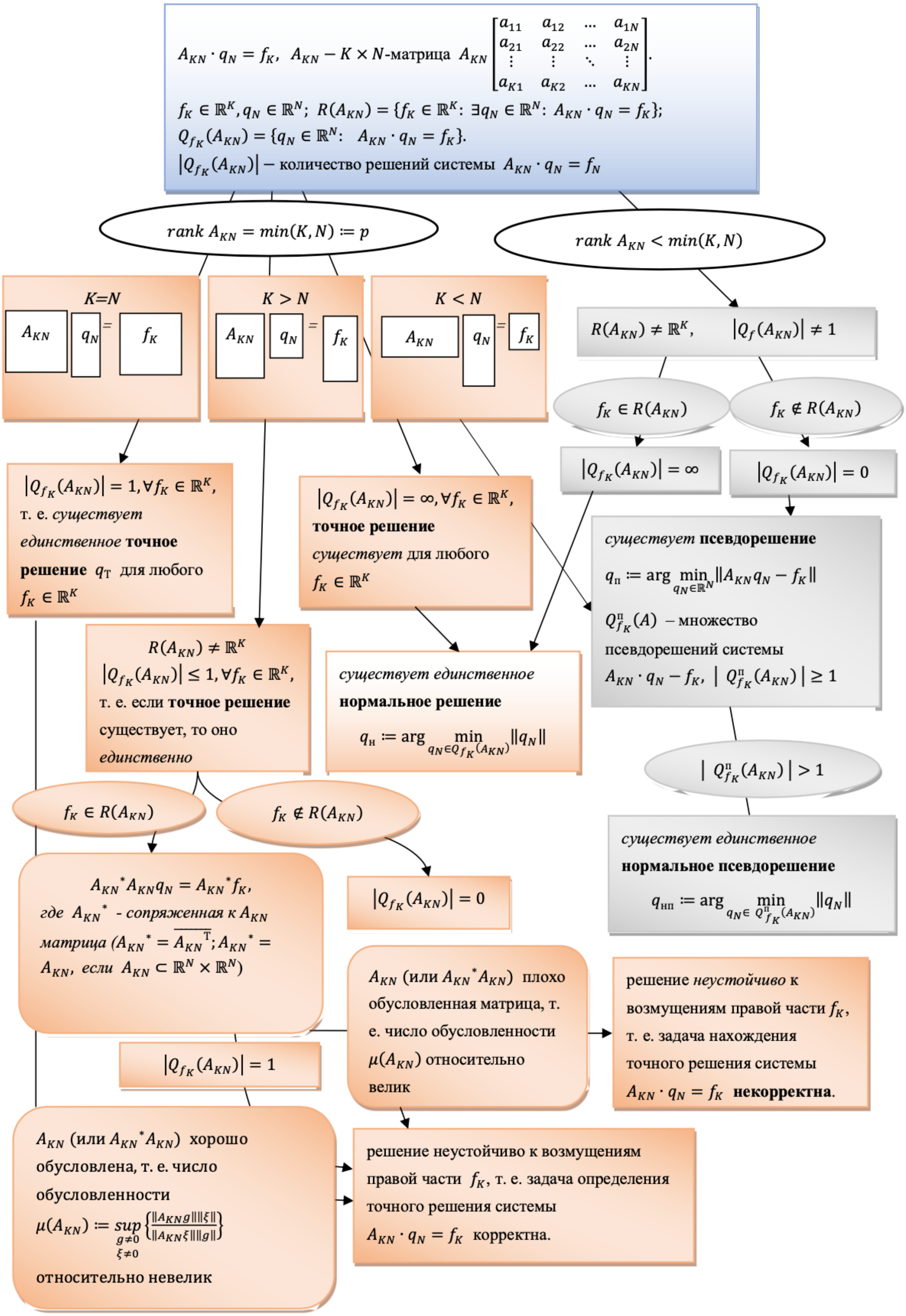

Рис. 11: Система линейных алгебраических уравнений (СЛАУ) (Кабанихин, 2023) [24].

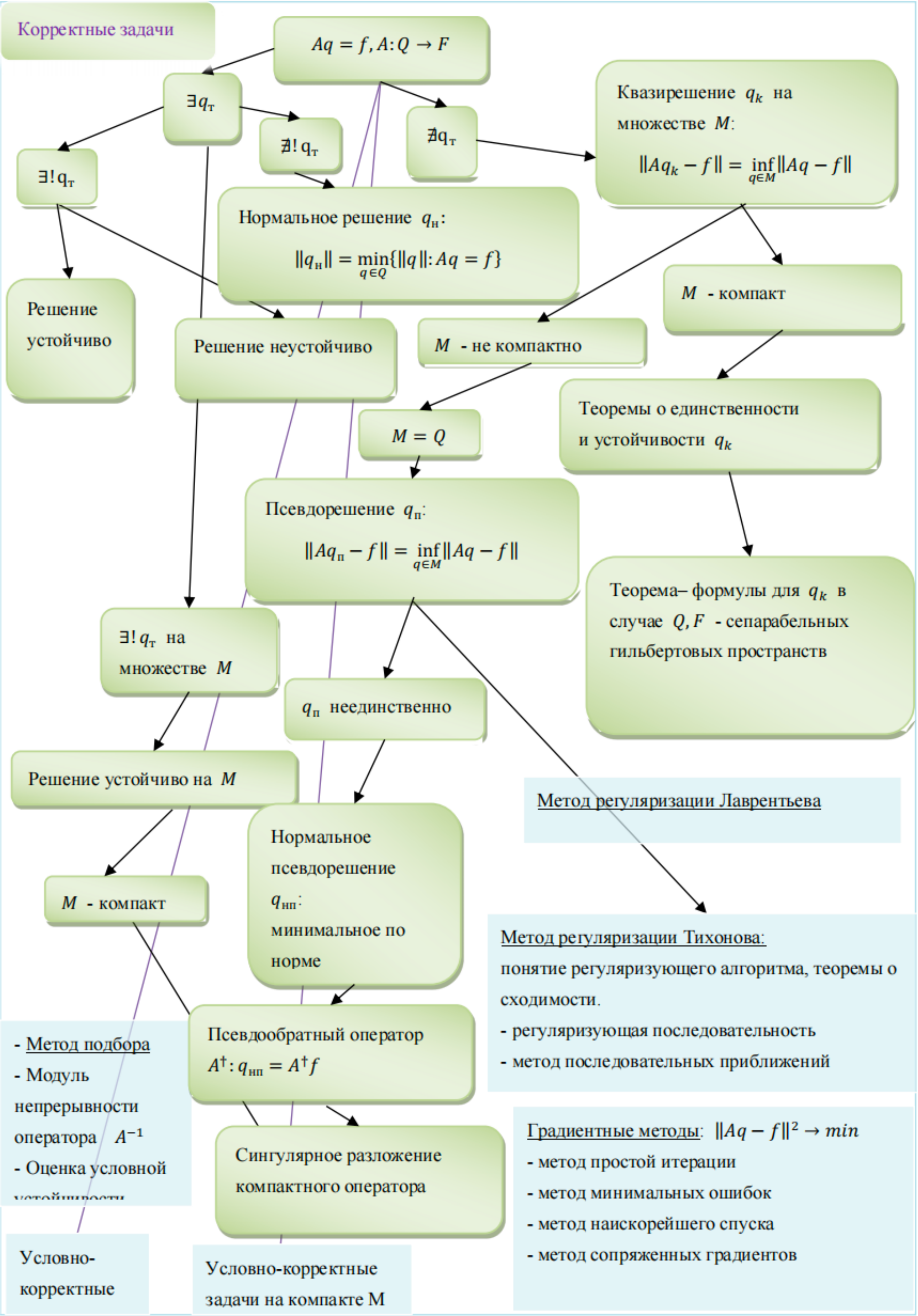

Рис. 12: Система корректных задач.

### 5.1. Псевдообращение матриц

Рассмотрим систему $K \times N$ линейных алгебраических уравнений с $n$ неизвестными

$$\begin{cases} a_{11}q_1 + a_{12}q_2 + \cdots + a_{1N}q_N = f_1, \\ a_{21}q_1 + a_{22}q_2 + \cdots + a_{2N}q_N = f_2 \\ \qquad \dots\dots\dots\dots\dots \\ a_{K1}q_1 + a_{K2}q_2 + \cdots + a_{KN}q_N = f_K \end{cases},$$

Эту систему можно записать в матричном виде

$$A_{KN} \cdot q_N = f_K, \tag{5.1}$$

где $A_{KN}$ — действительная $K \times N$-матрица коэффициентов системы, $f = (f_1, f_2, \cdots, f_K)^T \in \mathbb{R}^N$ — вектор-столбец правых частей, $q = (q_1, q_2, \cdots, q_N)^T \in \mathbb{R}^N$ — вектор-столбец неизвестных. Отметим, что в соответствии с принятыми в данной главе обозначениями $Q = \mathbb{R}^N$, $F = \mathbb{R}^K$, $A_{KN}: \mathbb{R}^N \longrightarrow \mathbb{R}^K$. Мы будем изучать методы построения приближенных решений задачи $A_{KN} \cdot q_N = f_K$ в случае, когда эта задача некорректна ($K \neq N$, или $A_{KN}$ вырождена, или $A_{KN}$ плохо обусловлена и т. п.).

Как хорошо известно, если $K = N$ и $\det A_{KN} \neq 0$, то существует единственное решение системы. В вырожденном случае ($\det A_{KN} = 0$) и в случае $K \neq N$ система может не иметь решения или иметь бесконечно много решений. Если система не имеет решения, то вводится понятие нормального псевдорешения так, что любая система с произвольной матрицей $A_{KN}$ будет обладать этим решением, причем единственным. В случае, когда система уравнений имеет много решений, среди решений выбирается минимальное по норме, которое называется *нормальным решением и обозначается* $q_н$.

Рассмотрим подробнее случай невырожденной квадратной матрицы $A_{KN}$. Как отмечалось выше, теоретически этот случай можно считать хорошим в смысле существования и единственности решения. Однако в теории вычислительных методов невырожденные матрицы подразделяют на две категории: “плохо обусловленные” и “хорошо обусловленные”. *Плохо обусловленными* называют матрицы, для которых решение системы уравнений практически является неустойчивым. Иначе говоря, небольшие погрешности в правой части системы или погрешности, неизбежно возникающие при численной реализации, приводят к существенному отклонению полученного решения от точного.

Одной из важных характеристик практической устойчивости решения системы линейных уравнений является число обусловленности. Числом обусловленности квадратной матрицы $A_{KN}$ называется величина

$$\mu(A_{KN}) = \sup_{\substack{q \neq 0 \\ \xi \neq 0}} \left\{ \frac{\|A_{KN}q\| \|\xi\|}{\|A\xi\| \|q\|} \right\}. \tag{5.2}$$

Для понимания влияния $\mu(A_{KN})$ на устойчивость, рассмотрим систему с возмущенной правой частью $A_{KN}(q_N + \delta q_N) = f_K + \delta f_K$. Здесь через $\delta q$ обозначено отклонение решения, вызванное возмущением $\delta f_K$ правой части. Очевидно, что $A_{KN}\delta q_N = \delta f_K$. Из определения числа обусловленности следует, что

$$\frac{\|\delta q_N\|}{\|q_N\|} \leq \mu(A_{KN}) \frac{\|A_{KN}\delta q_N\|}{\|A_{KN}q_N\|} = \mu(A_{KN}) \frac{\|\delta f_K\|}{\|f_K\|}. \tag{5.3}$$

и $\mu(A_{KN})$ — наименьшая константа, гарантирующая это неравенство. Таким образом, число обусловленности $\mu(A_{KN})$ позволяет оценить относительную погрешность решения $\|\delta q_N\|/\|q_N\|$ через относительную погрешность правой части $\|\delta f_K\|/\|f_K\|$. Плохо обусловленные системы — это системы с очень большим числом обусловленности $\mu(A_{KN})$.

Решение системы линейных алгебраических уравнений может быть некорректной задачей для прямоугольных матриц $A_{KN}$ $(K \neq N)$, а также для квадратных вырожденных или плохо обусловленных матриц.

Методы решения систем линейных уравнений делятся на две основные группы прямые и итерационные. *Прямые методы* дают решение после выполнения заранее известного числа операций. Эти методы сравнительно просты и универсальны, но требуют, как правило, большого объема памяти ЭВМ и накапливают погрешности в процессе решения, поскольку вычисления на любом этапе используют результаты предыдущих операций. В связи с этим прямые методы применимы для сравнительно небольших $(N < 200)$ систем с плотно заполненной матрицей и не близким к нулю определителем. К прямым методам относятся, например, метод Гаусса, метод Жордана, метод квадратного корня, метод сингулярного разложения и т. д. Прямые методы (иногда называемые точными, хотя ввиду неизбежных погрешностей при численной реализации это название является условным) также успешно применяются для решения систем линейных алгебраических уравнений.

Итерационные методы — это методы последовательных приближений. Объем вычислений заранее определить трудно, но они требуют меньшего объема памяти, чем прямые методы. Итерационные методы часто используются при регуляризации некорректных систем линейных уравнений. Однако, стоит отметить, что часто наиболее эффективным способом решения линейной системы является сочетание итерационного подхода, с прямыми методами. В таких смешанных алгоритмах итерационные методы используются для уточнения решений, полученных с помощью прямых методов.

### 5.2. Обобщение понятия решения. Псевдорешение

Прежде чем перейти к изложению методов приближенного решения системы уравнений $A_{KN} \cdot q_N = f_K$, обсудим понятие решения этой системы, которая в общем случае может быть переопределенной, переопределенной или плохо обусловленной.

Здесь и далее $q_N \in \mathbb{R}^N$, $f_K \in \mathbb{R}^K$, $A_{KN}$ — вещественная матрица размера $K \times N$, которую называем $K \times N$ - матрицей в случае, когда необходимо указать ее размер.

Вектор $q_\text{п} \in \mathbb{R}^N$, реализующий минимум нормы невязки

$$L(q_N) = \|A_{KN}q_N - f_K\|^2 \to min \qquad (5.4)$$

называется *псевдорешением* системы $A_{KN}q_N = f_K$,

т. е. $q_\text{п} = \arg \min\limits_{q_N \in \mathbb{R}^N} \|A_{KN}q_N - f_K\|^2$.

Поскольку для приращения функционала $L(q_N)$ справедливо представление

$$\delta L(q_N) = L(q_N + h) - L(q_N) = 2(A_{KN}q_N, A_{KN}h) - 2(A_{KN}h, f_K) + (A_{KN}h, A_{KN}h),$$

необходимым условием минимума будет

$$A_{KN}{}^T A_{KN}q_N - A_{KN}{}^T f_K = 0,$$

где $A_{KN}{}^T$ — транспонированная матрица. Следовательно, вектор $q_\text{п}$ являются решением системы уравнени

$$A_{KN}{}^T A_{KN}q_N = A_{KN}{}^T f_K. \qquad (5.5)$$

**Определение 5.1.** Система уравнений (5.5) называется *нормальной системой по отношению к системе $A_{KN} \cdot q_N = f_K$.*

Нетрудно доказать и обратное утверждение, а именно, что каждое решение системы (5.5) минимизирует невязку в (5.4) (Иванов, 1978) [24]. Таким образом, задачи (5.4) и (5.5) эквивалентны. Нетрудно убедиться, что задача (5.4) всегда имеет решение, хотя, возможно, не единственное. Поэтому в силу факта, эквивалентности, система (5.5) также разрешима для любых матриц $A_{KN}$ и векторов $f_K$. Таким образом, множество решений нормальной системы (5.5) совпадает с множеством псевдорешений системы $A_{KN} \cdot q_N = f_K$ . Обозначим это множество через $Q^{\text{п}}_{f_K}$.

Рассмотрим задачу отыскания точки минимума функционала

$$min\{\|q_N - q^0\|^2 : q_N \in Q^{\text{п}}_{f_K}\} \quad (5.6)$$

где $q^0$ — некоторый фиксированный вектор. Решение $q^0_{\text{нп}}$ задачи (5.6) существует и единственно, поскольку строго выпуклый функционал $\|\cdot\|^2_Q$ достигает на выпуклом замкнутом множестве $Q^{\text{п}}_{f_K}$; минимума в единственной точке.

**Определение 5.2.** Решение $q^0_{\text{нп}}$ задачи (5.6) будем называть нормальным относительно $q^0$ псевдорешением уравнения $A_{KN} \cdot q_N = f_K$. Нормальное относительно нулевого вектора, (наименьшее по норме) псевдорешение системы $A_{KN} \cdot q_N = f_K$называется нормальным псевдорешением этой системы (или нормальным обобщенным решением) и обозначается $q_{\text{нп}}$ : $q_{\text{нп}} = \arg \min_{q_{\text{п}} \in Q^{\text{п}}_{f_K}} \|q_{\text{п}}\|$.

Если система$A_{KN} \cdot q_N = f_K$разрешима, то нормальное относительно $q^0$ псевдорешение $q^0_{\text{нп}}$, совпадает с нормальным относительно $q^0$ решением этой системы, т. к. с решением, наименее уклоняющимся по норме от вектора $q^0$ .В частности, если система $A_{KN} \cdot q_N = f_K$ однозначно разрешима, то псевдорешение единственно и совпадает с обычным решением.

Нормальное псевдорешение существует, единственно и непрерывно зависит от ошибок в правой части $f_K$ , поскольку псевдообратный оператор, действующий в конечномерном пространстве, ограничен.

Нормальное псевдорешение $q_{\text{нп}}$ неустойчиво по отношению к возмущениям элементов матрицы. Неустойчивость нормального псевдорешения проиллюстрируем на следующем примере.

Предположим, нам необходимо решить СЛАУ вида $A_{KN} \cdot q_N = f_K$.

Если бы матрица $A_{KN}$ была квадратной и невырожденной (число уравнений равно числу неизвестных и все уравнения линейно независимы), то решение задавалось бы формулой $q_N = {A_{KN}}^{-1} f_K$.

Предположим, что число уравнений больше числа неизвестных, т. е. матрица $A_{KN}$ прямоугольная.

Домножим обе части уравнения на ${A_{KN}}^T$ слева,

$${A_{KN}}^T A_{KN} q_N = {A_{KN}}^T f_K \quad (5.7)$$

В левой части теперь квадратная матрица и ее можно перенести в правую часть.

Операция $q_N = \left({A_{KN}}^T A_{KN}\right)^{-1} {A_{KN}}^T f_K$ называется псевдообращением матрицы $A_{KN}$, $q_N$ — псевдорешением.

Нормальное псевдорешение:

Если матрица ${A_{KN}}^T A_{KN}$ вырождена, псевдорешений бесконечно много, причем найти их на компьютере нетривиально.

Для решения этой проблемы используется ридж-регуляризация матрицы ${A_{KN}}^T A_{KN}$

$${A_{KN}}^T A_{KN} + \lambda I. \quad (5.8)$$

где $I$ — единичная матрица, а $A_{KN}$ — коэффициент регуляризации. Такая матрица, невырождена для любых $\lambda > 0$.

Величина

$$q_N = \left( A_{KN}{}^T A_{KN} + \lambda I\right)^{-1} A_{KN}{}^T f_K. \quad (5.9)$$

называется нормальным псевдорешением. Оно всегда единственно и при небольших положительных $A_{KN}$ определяет псевдорешение с наименьшей нормой.

Псевдорешение соответствует точке, минимизирующей невязку, а нормальное псевдорешение отвечает псевдорешению с наименьшей нормой.

Заметим, что псевдообратная матрица $\left(A_{KN}{}^T A_{KN}\right)^{-1} A_{KN}{}^T$ совпадает с обратной матрицей $A_{KN}{}^{-1}$ в случае невырожденных квадратных матриц.

Дискретизация линейной обратной задачи обычно приводит к возникновению линейной системы уравнений

$$A_{KN} \cdot q_N = f_K,\ A_{KN} \in \mathbb{R}^{K\times N},\ q_N \in \mathbb{R}^{K\times N},\ f_K \in \mathbb{R}^{K\times N}. \quad (5.10)$$

с плохо обусловленной матрицей $A_{KN}$ плохо определенного ранга. Вычисление значимого приближенного решения линейной системы (5.9) в общем случае требует, чтобы система была заменена соседней системой, которая менее чувствительна к возмущениям. Эта замена называется регуляризацией. Регуляризация Тихонова - один из старейших и наиболее популярных методов регуляризации. В своей простейшей форме регуляризация Тихонова заменяет линейную систему на регуляризованную систему

$$\left(A_{KN}{}^T A_{KN} + \mu I\right) q_N = A^T f_K. \quad (5.11)$$

где $\mu \geq 0$ - параметр регуляризации, определяющий величину регуляризации, а $I$ - оператор идентификации. Для любого фиксированного $\mu > 0$ система (5.11) имеет единственное решение

$$q_\mu = (A^T A + \mu I)^{-1} A^T f_K. \quad (5.12)$$

Целью настоящей статьи является обсуждение нескольких итерационных методов для определения подходящего значения параметра регуляризации $\mu > 0$ и вычисления связанного с ним решения $q_\mu$ крупномасштабных задач вида (5.11). Отметим, что Бьорк описал, как итерационные методы для решения (5.11) могут быть модифицированы, чтобы быть применимыми к решению

$$\left(A_{KN}{}^T A_{KN} + \mu F^T F\right) q = A^T f_K.$$

для большого класса операторов регуляризации $F$.

Обратите внимание, что решение (5.12) из (5.11) удовлетворяет $q_\mu \to q_0 = A_{KN}{}^{\dagger} f_K$ как $\mu \searrow 0$ , где $A^{\dagger} f_K$ обозначает псевдоинверсию Мура–Пенроуза к $A_{KN}$. В интересующих нас задачах матрица $A_{KN}$ имеет много “крошечных” сингулярных значений, а вектор $f_K$ с правой стороны загрязнен ошибками измерения. Таким образом, решение Мура–Пенроуза в целом имеет “огромные” компоненты и представляет мало практического интереса.

Для дальнейшего использования отметим, что решение $g_\mu$ из (5.11) удовлетворяет задаче минимизации

$$\min_{q_N \in \mathbb{R}^{K\times N}} \{\|A_{KN} \cdot q_N - f_K\|^2 + \mu \|q_N\|^2\} \quad (5.13)$$

Псевдорешения в нейронных сетях означают нахождение приближенных решений задач оптимизации, которые могут быть достигнуты с помощью различных методов оптимизации при условии, что функция потерь имеет негладкую или невыпуклую форму. Псевдоразрешение может быть достигнуто с помощью различных методов оптимизации, таких как градиентный спуск, стохастический градиентный спуск и т. д. Однако псевдорешения не гарантируют, что будет найдено оптимальное решение, и могут привести к локальному минимальному значению функции потерь.

Численная линейная алгебра — это захватывающая область исследований, и большая их часть мотивирована проблемой, которая может быть сформулирована просто: дано решение, найти вектор решения такой, что $A_{KN} \cdot q_N = f_K$. Книги Варги и Янга содержат полные статьи об итерационных методах, использовавшихся в 1960-х и 1970-х годах. Книга Хаусхолдера дает достаточно хороший обзор итерационных методов - в частности, проекционных методов. К 1960 году методы последовательной перерелаксации (SOR) позволили эффективно решать системы в диапазоне 20 000 неизвестных, а к 1965 году системы порядка 100,000 можно было решать в задачах, связанных с вычислением собственных значений в кодах диффузии ядра. успех методов SOR привел к созданию богатой теории итерационных методов; она также могла быть плодотворно использована в более поздних работах. анализе методов (Saad, 2000)[25].

### 5.3. Метод градиентного спуска

Градиентный спуск: метод оптимизации, который используется для обучения линейных нейронных сетей. Он заключается в том, что на каждой итерации веса нейронов обновляются в направлении, противоположном градиенту функции потерь.

Одним из самых эффективных и широко применяемых методов регуляризации является итерационная регуляризация, основанная на минимизации функции потерь $L(q_N) = \|A_{KN} \cdot q_N - f_K\|^2$ и методах градиентного спуска. Основой градиентных методов является известное утверждение о том, что если в некоторой точке $q_N$ из $Q$ градиент функции $L(q_N)$ не равен нулю, то, двигаясь в направлении антиградиента в точку, можно уменьшить значение функционала $L(q_N)$ при условии, что шаг спуска достаточно мал. В самом деле, вспоминая определение градиента (полагаем для простоты, что $Q$ — гильбертово пространство), мы можем записать

$$L(q - \alpha L'q) - L(q) = \langle L'q, -\alpha L'q \rangle + o(\alpha\|L'q\|) = -\alpha\|L'q\|^2 + o(\alpha\|L'q\|). \qquad (5.14)$$

Здесь правая часть уравнения $\alpha$ достаточно мала. Таким образом, выбирая размер шаг, другим способом, либо корректируя направление спуска, либо строя не возрастающую последовательность минимизации в терминах функции. Более того, можно оценить скорость убывания функционала $L(q_n)$ даже в случае многих градиентных методов, когда задача $A_{KN} \cdot q_N = f_K$ некорректна. Случай сильной сходимости более сложен и, конечно, возможен только в том случае, если существует хотя бы одно точное решение уравнения $A_{KN} \cdot q_N = f_K$. Для последовательностей, построенных с помощью простых итераций, спуска и сопряженных градиентов, можно доказать $\|q_n - q_T\|$ монотонную сходимость к нулю. Если градиентный метод используется для решения нечетко определенной задачи, то последовательность минимизации сначала приближается к точному решению, а затем, по мере увеличения числа итераций, значения могут начать расти. Стоит отметить, как выбираются параметры $n$ и где останавливаются итерации. Наиболее часто на практике используется принцип неограниченности, основанный на естественном предположении, что если неограниченность достигла уровня неопределенности измерения $L$, то нет необходимости продолжать. То же самое относится и к выбору параметров регуляризации $\alpha$. Однако, если удастся более глубоко исследовать задачу $A_{KN} \cdot q_N = f_K$ и получить оценку условной устойчивости, можно, во-первых, оценить сильную скорость сходимости, а во-вторых, сформулировать новое правило для выбора числа остановки итерационного процесса (Тихонов, 1943) [26].

Алгоритмы обучения линейных нейронных сетей обычно включают следующие: метод градиентного спуска (Fedosin et al., 2010) [27], метод обратного распространения ошибки (Romanov, 2009) [28], метод стохастического градиентного спуска (Bottou, 2012) [29], метод адаптивного градиента (Duchi et al., 2011) [20].

### 5.4. Итерационные регуляризирующие алгоритмы

Для решения плохо обусловленных систем $A_{KN} \cdot q_N = f_K$ могут быть использованы итерационные методы. Роль параметра регуляризации играет число итераций. Исследуем этот вопрос на примере метода простой итерации

$$q_{n+1} = (I - A_{KN}{}^T A_{KN}) q_n + A_{KN}{}^T A_{KN}, \; n = 0,1,2 \ldots N. \qquad (5.15)$$

В случае приближенных данных $\{A_h, f_\delta\}$ вместо (5.15) имеем итерационную схему.

$$q_{n+1} = (I - A_h^T A_h) q_n + A_h^T f_\delta, \; n = 0,1,2 \ldots N, \qquad (5.16)$$

где $\|A - A_h\| \leq h$, $\|f - f_\delta\| \leq \delta$.

Будем предполагать, что система $A_{KN} \cdot q_N = f_K$ разрешима, $\|A_{KN}\| \leq 1$, $\|A_h\| \leq 1$. Последнее условие не ограничивает класса решаемых задач, поскольку его выполнения можно всегда достичь, умножив уравнение $Aq = f$ на подходящую константу. Как и ранее, через $q_{\text{н}}^0$ обозначаем нормальное относительно $q_0$ решение, т. е. решение системы $A_{KN} \cdot q_N = f_K$, для которого норма $\|q_H^0 - q_0\|$ минимальна, где $q_H^0$ начальное приближение в процессах (5.15), (5.16).

Определим правила выбора параметра $n$, при котором следует остановить процесс (5.13):

Правило 1. Зададим числа $a_1 > 0$, $a_2 > 0$ и выбираем такое $n(h, \delta)$, для которого впервые выполнено $\|q_n - q_{n-1}\| \leq a_1 h + a_2 \delta$.

Правило 2. Зададим числа $a_0 \geq \|q_H^0\|$, $a_1 > 1$ и выберем такой номер $n(h, \delta)$, для которого впервые выполнено $\|A_h q_n - f_\delta\| \leq a_0 h + a_1 \delta$. (Разумеется, в данном случае предполагается, что нам известна верхняя оценка для $\|q_H^0 - q_0\|$!)

Правило 3. Зададим числа $a_1 > 1$, $a_2 > 1$ и $a > 1$ и выберем такой номер $n$, для которого впервые будет выполнено хотя бы одно из неравенств:

$$\|A_h q_n - f_\delta\| \leq a_1 \|q_n\| h + a_2 \delta, \; n \geq a / (a_1 \|q_n\| h + a_2 \delta)^2.$$

Имеет место следующая теорема (Вайникко, 1980) [32].

**Теорема 5.1.** *Пусть последовательные приближения* (5.16) *останавливаются по любому из правил* 1,2 *или* 3. Тогда

$$\lim_{h,\delta \to 0} \|q_{n(h,\delta)} - q_{\text{н}}^0\| = 0. \qquad (5.17)$$

*При этом для числа итераций* $n(h, \delta)$ *в случае правила* 1 *справедливо Соотношение*

*$(h + \delta) n(h, \delta) \longrightarrow 0$ при $h, \delta \longrightarrow 0$,*

*а в случае правил* 2 *и* 3 — *соотношение*

*$(h + \delta)^2 n(h, \delta) \longrightarrow 0$ при $h, \delta \longrightarrow 0$.*

Соотношение (5.14) означает, что правила 1, 2 и 3 определяют регуляризирующие алгоритмы решения уравнения $A_{KN} \cdot q_N = f_K$.

### 5.5. Обоснование метода регуляризации

Рассмотрим понятие регуляризирующего алгоритма и связанного с ним понятия регуляризованного семейства приближенных решений, введенного А. Н. Тихоновым (Тихонов,1963) [33]. Всюду в данном разделе $A_{KN}$ — квадратная матрица. Условимся различать в задаче $A_{KN} \cdot q_N = f_K$ точные данные $\{A_{KN}, f_K\}$ , которые нам неизвестны, и приближенные данные $\{A_h, f_\delta\}$ , $h > 0$ , $\delta > 0$ , $\|A_h - A\| \le h$ , $\|f_\delta - f_K\| \le \delta$ с уровнем погрешности $h$ , $\delta$ . Рассмотрим способы построения семейства векторов, сходящихся к нормальному псевдорешению $q_{\text{нп}}$ уравнения $A_{KN} \cdot q_N = f_K$ при $h, \delta \to 0$. За приближенное решение нельзя, вообще говоря, принимать нормальное псевдорешение возмущенного уравнения

$$A_{hq} = A_{KN}{}^T f_K \qquad (5.18)$$

поскольку нормальное псевдорешение неустойчиво.

Пусть в нашем распоряжении имеется семейство алгоритмов $\{R_\alpha\}_{\alpha>0}$, каждый из которых паре $\{A_h, f_\delta\}$, $h > 0$, $\delta > 0$, однозначно сопоставляет вектор $R_\alpha(A_h, f_\delta) = q_{h\delta}^\alpha \in Q = R^N$ . Если существует зависимость параметра $\alpha = \alpha(h, \delta)$ от погрешностей $\delta, h$ исходных данных, такая, что $\lim\limits_{h,\delta\to 0} \alpha(h, \delta) = 0$,

$$\lim_{h,\delta\to 0} \|q_{h\delta}^\alpha - q_{\text{нп}}\| = 0. \qquad (5.19)$$

То $\{q_{h\delta}^\alpha\}$ , $h > 0$ , $\alpha > 0$ , $\delta > 0$ , называется *регуляризованным семейством* приближенных решений, а алгоритм $R_\alpha$ — *регуляризирующим алгоритмом* для задачи $A_{KN} \cdot q_N = f_K$. Напомним, что нормальное псевдорешение $q_{\text{нп}}$ системы $A_{KN} \cdot q_N = f_K$ совпадает с ее точным решением, если эта система однозначно разрешима, и с нормальным решением, если система $A_{KN} \cdot q_N = f_K$ имеет множество решений.

Если $A_{KN}{}^{-1}$ существует, т. е. $A_{KN}$ — невырожденная матрица, то для достаточно малых $h$ обратная матрица $A_h^{-1}$тоже существует и решение уравнения (5.18) теоретически будет сходиться к решению уравнения $A_{KN} \cdot q_N = f_K$при $h, \delta \to 0$. Однако, если $A_{KN}$— плохо обусловленная матрица, т. е. ее число обусловленности $\mu(A_{KN}) = \sup\limits_{\substack{g\ne 0\\ \xi\ne 0}} \left\{\frac{\|A_{KN} g\| \|\xi\|}{\|A\xi\| \|g\|}\right\}$ велико, то уклонение решения возмущенной системы (5.18) от нормального псевдорешения системы $A_{KN} \cdot q_N = f_K$ даже при малых $h, \delta$ может оказаться недопустимо большим и задачу следует считать практически неустойчивой (некорректной).

Перейдем к описанию конкретных процедур построения регуляризованных приближенных решений $q_{h\delta}^\alpha$.

Рассмотрим сначала частный случай, когда $A_{KN}$ — симметричная положительно полуопределенная матрица, для которой система $A_{KN} \cdot q_N = f_K$ при данном векторе $f_K$ разрешима. Применим к этому случаю схему регуляризации М. М. Лаврентьева.

Перейдем к регуляризованной системе

$$(A_{KN} + \alpha I)q = f_K + \alpha q^0. \qquad (5.20)$$

где $\alpha$ — положительный параметр, $I$ — единичная матрица, $q^0$ — пробное решение, т. е. некоторое приближение к искомому решению (если информации о решении нет, то можно положить $q^0 = 0$).

При принятых условиях система (5.20) имеет единственное решение $q^\alpha$, которое сходится при $\alpha \to 0$ к нормальному относительно $q^0$решению $q_H^0$. Справедлива следующая лемма 4.1.

**Лемма 5.1.** *Пусть* $\{A_h, f_\delta\}$, $h > 0$, $\delta > 0$, *таковы, что* $\|A_h - A_{KN}\| \le h$, $\|f_\delta - f_K\| \le \delta$, *и* $A_h$ *— симметричная положительно полуопределенная матрица. Тогда система*

$$(A_{KN} + \alpha I)q = f_\delta + \alpha q^0. \tag{5.21}$$

*однозначно разрешима и, если* $\alpha(h,\delta) \to 0$, $(h,\delta)/\alpha(h,\delta) \to 0$ *при* $h,\delta \to 0$, *то ее решение* $q_{h\delta}^\alpha$*; сходится к нормальному относительно* $q^0$ *решению* $q_H^0$ *уравнения* $A_{KN} \cdot q_N = f_K$, *т. е. к решению, наименее уклоняющемуся от вектора* $q^0$.

Согласно решение систем (5.21) $\{q^{\alpha(h,\delta)}\}$ образуют регуляризованное семейство приближенных решений для системы $A_{KN} \cdot q_N = f_K$, причем выбор параметра по формуле $\alpha = \sqrt[p]{h+\delta}$ $(p > 1)$удовлетворяет необходимым требованиям, поскольку $\alpha = \sqrt[p]{h+\delta} \to 0$, $(h,\delta)/\sqrt[p]{h+\delta} = (h+\delta)^{1-1/p} \to 0$ при $h,\delta \to 0$.

**Лемма 5.2.** *Пусть* $\|(A_{KN} + \alpha B)^{-1} A_{KN}\| \le C < \infty$ (*при* $\alpha \longrightarrow 0$ ), $\|A_h - A_{KN}\| \le h$, $\|B_\mu - B\| \le \mu$, $\|f_\delta - f_K\| \le \delta \|f_K\|$. *Тогда при достаточно малых h, δ система*

$$(A_h + \alpha B_\mu)q = f_\delta. \tag{5.22}$$

имеет единственное решение $q_{h\delta\mu}^\alpha$ и справедлива оценка погрешности

$$\|q_{h\delta\mu}^\alpha - q_H^B\| \le C\left(|\alpha| + \mu + \frac{h+\delta}{|\alpha|}\right). \tag{5.23}$$

где $q_H^B$ — решение системы $A_{KN} \cdot q_N = f_K$, удовлетворяющее условию

$$\|Bq_H^B\| = \min\{\|Bq : q \in Q_f\|\}.$$

($Q_{f_K}$ — множество решений системы $A_{KN} \cdot q_N = f_K$).

Из оценки (5.23) вытекает, что если $\alpha(h,\delta,\mu) \longrightarrow 0$ и $(h+\delta)/|\alpha(h,\delta,\mu)| \longrightarrow 0$ при $h,\delta,\mu \longrightarrow 0$, то имеет место сходимость

$$\lim_{h,\delta,\mu \to 0} \|q_{h\delta\mu}^\alpha - q_H^B\| = 0.$$

Наиболее важным моментом в описанной регуляризации является подбор матрицы $B$, для которой $A_{KN}$+$\alpha_0 B$ невырождена и $\|( A_{KN} + \alpha B)^{-1} A_{KN}\| < \infty$. В работе А. Б. Назимова (Назимов и Морозов, 1985) [34] можно найти способы построения матриц с таким свойством.

Исследуем, наконец, общую ситуацию, когда система $A_{KN} \cdot q_N = f_K$, вообще говоря, неразрешима. В этом случае искомым является нормальное относительно некоторого $q^0$ псевдорешение $q_{\text{нп}}^0$. Решается задача устойчивой аппроксимации $q_{\text{нп}}^0$, когда $A_{KN}$ и $f_K$ заданы с погрешностью. Нормальное относительно $q^0$ псевдорешение $q_{\text{нп}}^0$ неустойчиво к возмущениям элементов матрицы, поэтому необходимо использовать регуляризацию. В качестве регуляризованного приближенного решения выберем вектор $q_{h\delta}^\alpha$, удовлетворяющий системе

$$(A_h^T A_h + \alpha I)q_N = A_h^T f_\delta + \alpha q^0 \tag{5.24}$$

**Лемма 5.3.** *Пусть* $\|A_h - A_{KN}\| \le h$, $\|f_\delta - f_K\| \le \delta$, $\alpha > 0$. *Тогда система (5.24) однозначно разрешима и справедлива оценка*

$$\|q^0_{\text{нп}} - q^\alpha_{h\delta}\| \le c_1\alpha + \tfrac{h}{\alpha}(\|Aq^0_{\text{нп}} - f_K\| + 2c_2^2\alpha^2)^{1/2} + \tfrac{1}{\alpha^{1/2}}(c_3 h + \delta). \quad (5.25)$$

*где* $c_i$ $(i = 1,2,3)$ *— константы, которые зависят от нормы* $\|q^0_{\text{нп}}\|$ *нормального относительно* $q^0$ *псевдорешения.*

**Следствие 5.1.** *Пусть* $h$, $\delta$ *суть величины порядка* $\varepsilon$, *причем* $\varepsilon$ *достаточно мало. Если уравнение* $A_{KN} \cdot q_N = f_K$ *имеет точное решение (т. е.* $\|Aq^0_{\text{нп}} - f_K\| = 0$), *то правая часть оценки* (5.25) *по характеру зависимости от* $\alpha$ *и* $\varepsilon$ *есть функция вида*

$$\varphi(\alpha,\varepsilon) = \alpha + \varepsilon + \frac{\varepsilon}{\alpha^{1/2}}. \quad (5.26)$$

*При* $\alpha = \varepsilon^{2/3}$ *она принимает значение порядка* $\varepsilon^{2/3}$. *Если система* $f_K$ *не разрешима* ($\|Aq^0_{\text{нп}} - f_K\| \ne 0$), *то правая часть неравенства* (5.25) *есть функция вида*

$$\psi(\alpha,\varepsilon) = \alpha + \frac{\varepsilon}{\alpha} + \frac{\varepsilon}{\alpha^{1/2}}. \quad (5.27)$$

*При* $\alpha = \varepsilon^{1/2}$ *она принимает значение порядка* $\varepsilon^{1/2}$. *Эти оценки получаются минимизацией функций* $\varphi$ *и* $\psi$ *по* $\alpha$ *(т. е. решением уравнений* $\varphi'_\alpha(g) = 0, \psi'_\alpha(g) = 0$).

Таким образом, если входные данные системы $A_{KN} \cdot q_N = f_K$ заданы с точностью порядка $\varepsilon$, то нормальное относительно $q^0$ решение может быть определено с точностью порядка $\varepsilon^{2/3}$ в случае разрешимости точного уравнения $A_{KN} \cdot q_N = f_K$; в противном случае можно построить нормальное относительно $q^0$ псевдорешение с точностью порядка $\varepsilon^{1/2}$.

Заметим, что более сложный способ выбора параметра $\alpha$ позволяет аппроксимировать нормальное относительно $q^0$ псевдорешение с точностью порядка аппроксимации $h + \delta$ (Джумаев,1982) [35].

Задача (5.24) эквивалентна, задаче на минимум

$$\min\left\{\left\|A_{hq} - f_\delta\right\|^2 + \alpha\|q_N - q^0\|^2 : q_N \in \mathbb{R}^n\right\}. \quad (5.28)$$

В такой форме способ построения $q^\alpha_{h\delta}$ ; известен как вариационный метод регуляризации А. Н. Тихонова. В некоторых случаях целесообразно использовать более общий вид стабилизирующего функционала, а именно

$$\min\left\{\left\|A_{hq} - f_\delta\right\|^2 + \alpha\|L(q_N - q^0)\|^2 : q_N \in \mathbb{R}^N\right\}, \quad (5.29)$$

где $L$ — невырожденная квадратная матрица, подходящий выбор которой позволяет повысить точность регуляризованного решения.

### 5.6. Принципы выбора параметра регуляризации

Параметр регуляризации $\alpha$ необходимо связывать с погрешностями $\delta$, $h$ исходных данных: $\alpha = \alpha(\delta, h)$. В прикладных задачах обычно уровень погрешности $\delta$, $h$ фиксирован ($\delta$, $h$ не стремятся к нулю) и нужно указать конкретное $\alpha(\delta, h)$, в определеном смысле, наилучшее. Дело в том, что при уменьшении $\alpha$ ухудшается обусловленность матриц регуляризованных систем и, следовательно, могут возникнуть вычислительные погрешности, а при увеличении $\alpha$ приближенное решение плохо аппроксимирует точное решение. Поэтому здесь необходим разумный компромисс.

Опишем некоторые принципы выбора параметра $\alpha$.

Принцип невязки ($\alpha_н$). Предположим, что погрешность имеется лишь в правой части системы $A_{KN} \cdot q_N = f_K$ , т. е. $A_h \equiv A_{KN}(h = 0)$ , $\|f_\delta - f_K\| \le \delta$ . Обозначим через $A_{KN}q_\delta^\alpha$ решение системы при $h = 0, q^0 = 0$. Значение параметра $\alpha$ в принципе невязки выбирается таким образом, чтобы выполнялось соотношение (Тихонов, 1974) [36] (Морозов, 1987) [37].

$$\|A_{KN}q_\delta^\alpha - f_\delta\| = \delta. \qquad (5.30)$$

Принцип обобщенной невязки ($\alpha_{он}$). Данный принцип охватывает общий случай задания погрешности: $\|A_{KN} - A_h\| \le h$, $\|f_K - f_\delta\| \le \delta$. Параметр $\alpha$ находится из уравнения

$$\|A_h q_{h\delta}^\alpha - f_\delta\| = h\|q_{h\delta}^\alpha\| + \delta, \qquad (5.31)$$

для решения которого разработаны численные методы и программы (Тихонов,1983) [38].

**5.7. Численные эксперименты**

**Пример 5.7.1.** Рассмотрим несовместную систему

$$\begin{cases} 1q_1 + 0q_2 = 1, \\ 0q_1 + 0q_2 = 1, \end{cases} A = \begin{bmatrix} 1 & 0 \\ 0 & 0 \end{bmatrix}, f = \begin{bmatrix} 1 \\ 1 \end{bmatrix}. \qquad (5.32)$$

Перейдем к нормальной системе с

$$A^T A = \begin{bmatrix} 1 & 0 \\ 0 & 0 \end{bmatrix}, A^T f = \begin{bmatrix} 1 \\ 0 \end{bmatrix}. \qquad (5.33)$$

и построим множество псевдорешений системы (5.32):

$$Q_f^п = \{(q_1, q_2): q_1 = 1, q_2 \in \mathbb{R}\}. \qquad (5.34)$$

Тогда нормальным псевдорешением системы $Aq = f$ будет вектор $q_{нп} = (1{,}0)^T$ .

Возмущенную систему возьмем в виде

$$\begin{cases} 1q_1 + 0q_2 = 1, \\ 0q_1 + hq_2 = 1, \end{cases} \qquad (5.35)$$

где $h$ — малый параметр. Поскольку эта. система имеет единственное решение $q_{нпh} = (1{,}1/h)^T$, оно и будет единственным (и, значит, нормальным) псевдорешением возмущенной системы. Очевидно, $q_{нпh} \to \infty$ при $h \to 0$, т. е, что доказывает неустойчивость нормального псевдорешения по отношению к возмущениям элементов матрицы.

**Пример 5.7.2.** Рассмотрим несовместную систему

$$\begin{cases} 1q_1 + 0q_2 + 0q_3 = 1, \\ 0q_1 + 1q_2 + 0q_3 = 1, \\ 0q_1 + 0q_2 + 0q_3 = 1. \end{cases} A = \begin{bmatrix} 1 & 0 & 0 \\ 0 & 1 & 0 \\ 0 & 0 & 0 \end{bmatrix}, f = \begin{bmatrix} 1 \\ 1 \\ 1 \end{bmatrix}. \qquad (5.36)$$

В следующем примере используется метода SVD для поиска псевдорешения системы линейных уравнений (5.36).

Сначала мы определяем матрицу коэффициентов $A$ и результирующий вектор $f$. Затем мы используем функцию SVD numpy для вычисления SVD-разложения $A$. Затем мы вычисляем псевдоинверсную матрицу на основе SVD-разложения и умножаем ее на результирующий вектор $f$, чтобы получить псевдорешение $q_{нп} = (1, 1, 0)^T$ системы линейных уравнений (5.36).

**Пример 5.7.3.** Рассмотрим несовместную систему

$$\begin{cases} 1q_1 + 0q_2 + 0q_3 = 1, \\ 0q_1 + 10q_2 + 0q_3 = 1, \\ 0q_1 + 0q_2 + 0q_3 = 1. \end{cases} A = \begin{bmatrix} 1 & 0 & 0 \\ 0 & 10 & 0 \\ 0 & 0 & 0 \end{bmatrix}, f = \begin{bmatrix} 1 \\ 1 \\ 1 \end{bmatrix}. \qquad (5.37)$$

В следующем примере используется 20-слойных линейных нейронных сетей для поиска псевдорешения системы линейных уравнений (5.37). Тогда нормальным псевдорешением уравнений (5.37) будет вектор $q_{нп} = (1, 0.1, 0)^T$ .

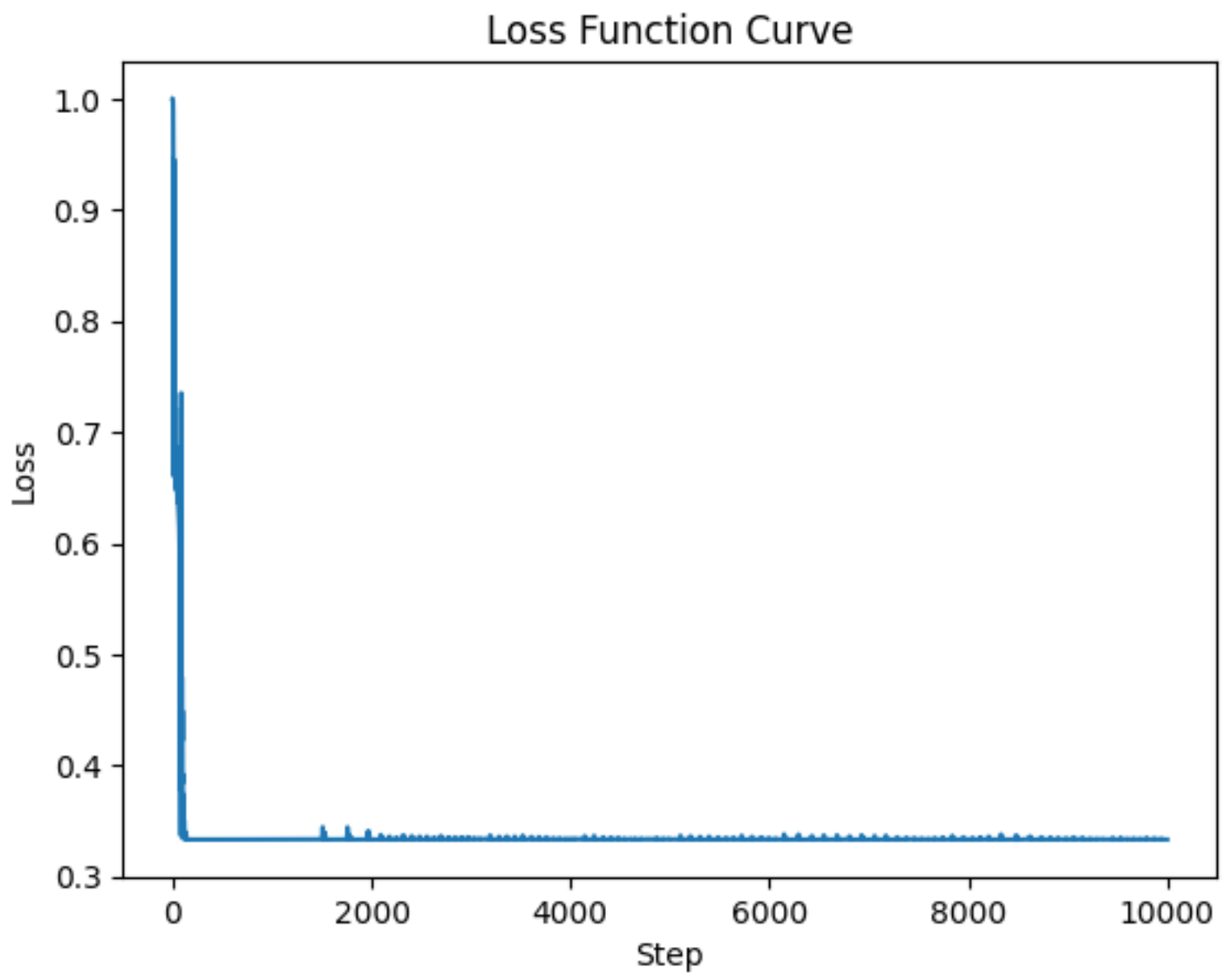

Рис. 13: Вывод значения функции потерь на каждой итерации обучения.

## 6. Заключение

В отличие от нелинейных нейронных сетей, линейные сети не могут решать сложные задачи, такие как распознавание образов или классификация изображений. Линейные нейронные сети подходят для задач регрессии, то есть люди хотят предсказать значение на основе входных данных. Линейные нейронные сети также могут использоваться в качестве базовых моделей для более сложных архитектур, таких как сверточные нейронные сети и рекуррентные нейронные сети.

Регуляризация также используется в нейронных сетях для решения проблемы чрезмерного веса сети и чрезмерного обучения. Использование L1-регуляризации не имеет значения "выбора признаков", как в линейных моделях. Регуляризация является неотъемлемой частью глубокого обучения нейронных сетей. Они либо ограничивают параметры обучения, либо вносят помехи где-то в цикле обучения, в зависимости от данных обучения, архитектуры сети, параметров обучения или целей. Однако регуляризация не уменьшает количество параметров и не упрощает структуру сети. Для нейронных сетей, помимо добавления штрафных санкций к эмпирическим рискам, также активно применяется другой метод борьбы с переобучением — отсев (dropout). В этом процессе сеть будет упрощена в соответствии с правилами - если функция ошибки не изменится, сеть может быть упрощена еще больше. Кроме того, стандартизация партии и ранняя остановка также могут решить проблему исчезновения градиента /взрыва.

**Список литературы**

## 7. Приложение

### 7.1. Краткий исторический обзор

Исследования по линейной регрессии появились в начале XIX века (Legendre, 1805; Gauss, 1809) [39-40]. Связь между независимой переменной $x$ и зависимой переменной $y$ предполагалась линейной.

В 1911 году Герман Вейль опубликовал один из первых примеров решения обратных задач, описывающий асимптотическое поведение собственных значений оператора Лапласа-Бельтрами (Weyl, 1911) [41]. В 1936 году Алан Тьюринг предложил концепцию машины Тьюринга, которая положила начало эре искусственного интеллекта (Turing, 1936) [42]. В 1983 году Юрий Нестеров предложил версию градиентного спуска, которая сходилась намного быстрее, чем обычный метод градиентного спуска (Nesterov,1983) [43]. В 1997 году Фукумидзу Кэндзи исследовал линейные нейронные сети (Fukumizu, 1997) [2]. В 2012 году Джон Каливас описал варианты регуляризации L2 и L1, А. Н. Тихонова (Kalivas, 2012) [9]. В 2016 году Ян Гудфеллоу, Йошуа Бенгио и Аарон Курвиль в своей книге "Глубокое обучение" отметили: "В контексте глубокого обучения большинство стратегий регуляризации основаны на регуляризирующих оценках. Регуляризация оценки достигается за счет уменьшения дисперсии с увеличением смещения. Эффективный регуляризатор может значительно уменьшить дисперсию без чрезмерного увеличения смещения". (Goodfellow, 2016) [44].

В 2021 году Хаджи Саад и Абдулазиз Мохсинз сравнили различные методы оптимизации, основанные на алгоритмах градиентного спуска (Haji and Abdulazeez, 2012) [45].

Основателями и теоретиками нейронных сетей считаются Джеффри Хинтон (Rumelhart, 1985) [16], Ян Лекун (LeCun, 1995) [19] и Джошуа Бенгио (Bengio, 2000) [46]. Их работы стали основой для многих современных разработок в области технологий искусственного интеллекта, и они по-прежнему остаются известными и влиятельными учеными в данной области.

Джеффри Хинтон впервые предложил алгоритм "обратного распространения ошибок", который помогает обучать нейросети в многослойных сложных сценариях, а также знаменитую технологию "отсева", которая может предотвратить проблему трансформации модели (Rumelhart and McClelland,1987; McClelland and Rumelhart, 1987) [47-48].

### 7.2. Теорема о сингулярном разложении

Пусть задана система линейных алгебраических уравнений (Годунов,1966) [49] (Годунов,1980) [50].

$$A_{KN} \cdot q_N = f_N \tag{7.1}$$

где $A_{KN}$ — действительная $K \times N$ -матрица $q_N \in \mathbb{R}^N$ , $f_N \in \mathbb{R}^K$ , и требуется решить эту систему относительно неизвестного вектора $q$ , любое ортогональное преобразование задается ортогональной матрицей, т. е. такой матрицей $U$, что $U^T U = U U^T = I$, где $U^T$ — транспонированная матрица, $I$ — единичная матрица. Свойство сохранения нормы вектора при ортогональных преобразованиях позволяет искать псведорешения вырожденных систем путем замены исходной задачи минимизации невязки $\|A_{KN} \cdot q_N - f_N\|$ задачей минимизации функционала $\|U^T(Aq - f)\|$ , в которой матрица $U^T A$ имеет более простую (например, блочную) структуру благодаря специальному построению ортогональной матрицы $U$.

Самым известным из полных ортогональных разложений является *сингулярное разложение* $K \times N$-матрицы $A$, т. е. разложение вида,

$$A_{KN} = U_{KK}\Sigma_{KN} V_{NN}^{*} \tag{7.2}$$

где $V$ — ортогональная $N \times N$ -матрица, ( $U$ — ортогональная $K \times K$ -матрица, $\Sigma$ — диагональная $K \times N$ -матрица, у которой $\sigma_{ij} = 0$ при $i \neq j$ ; $\sigma_{ij} = \sigma_i \geq 0$ . Величины $\sigma_i$ ; называются сингулярными числами матрицы $A$. Всюду считаем, что $\sigma_i$; занумерованы в порядке невозрастания $\sigma_{i+1} \leq \sigma_i$.

Напомним, что *диагональной* $K \times N$ *—матрицей* $\Sigma_{KN}$ называют $K \times N$—матрицу с элементами $\sigma_{ij}$, удовлетворяющими условиям

$$\sigma_{ij} = \begin{cases} 0, i \neq j, \\ \sigma_i, i = j. \end{cases}$$

Например, при $N > K$ диагональная $K \times N$-матрица $\Sigma_{KN}$ имеет вид

$$\Sigma_{KN} = \begin{bmatrix} \sigma_1 & 0 & 0 & \dots & 0 \\ 0 & \sigma_2 & 0 & & 0 \\ & \vdots & & \ddots & \vdots \\ 0 & 0 & \sigma_\rho & \cdots & 0 \end{bmatrix}$$

И в том, и в другом случае диагональную $K \times N$—матрицу для краткости будем обозначать $\Sigma_{KN}$ =diag$(\sigma_1, \sigma_2, \cdots, \sigma_\rho)$, $\rho = min\{K, N\}$ (Кабанихин, 2021) [51].

**Теорема 7.1. (о сингулярном разложении).** *Для любой вещественной* $K \times N$ *— матрицы* $A_{KN}$ *можно подобрать ортогональные* $K \times K$*—матрицу* $U$ *и* $N \times N$*—матрицу* $V$, *а также диагональную* $\Sigma_{KN}$ *—матрицу*

$$\Sigma_{KN} = \text{diag}(\sigma_1, \sigma_2, \cdots, \sigma_\rho), \rho = min\{K, N\}.$$

такие, что

$$A_{KN} = U_{KK}\Sigma_{KN} V_{NN}^{*}$$

$$0 \leq \sigma_\rho \leq \sigma_{\rho-1} \leq \cdots \leq \sigma_2 \leq \sigma_1.$$

Числа $\sigma_i = \sigma_2(A_{KN})$, $i = \overline{1, \rho}$, определяются однозначно и называются сингулярными числами матрицы $A$.

**Лемма 7.1. (о взаимной ортогональности сингулярных векторов)**. *Пусть* $u_{(k)}$ *и* $v_{(k)}$ *— соответственно правые и левые сингулярные векторы матрицы* $A_{KN}$ , *тогда* $\langle u_{(k)}, u_{(j)} \rangle = 0$ *и* $\langle u_{(k)}, v_{(j)} \rangle = 0$ *при* $k \neq j$.

$$A_{KN}{}^T A_{KN} v_{(k)} = A_{KN}{}^T \sigma_{(k)} u_{(k)} = \sigma_{(k)}^2 v_{(k)}, \; A_{KN} A_{KN}{}^T u_{(k)} = \sigma_{(k)}^2 u_{(k)},$$

*и, следовательно, правые сингулярные векторы* $K \times N$ *—матрицы* $A_{KN}$ *являются собственными векторами* $N \times N$ *—матрицы* $A_{KN}{}^T A_{KN}$ , *а левые —* $K \times N$ *—матрицы* $A_{KN} A_{KN}{}^T$.

**Теорема 7.2.** *Для каждой* $K \times N$*—матрицы* $A_{KN}$ *существует ортонормированная систсма из* $n$ *правых сингулярных векторов и ортонормированная система из* $K$ *левых сингулярных векторов, которые называются сингулярными базисами матрицы* $A_{KN}$.

**Лемма 7.4.** *Квадратная* $N \times N$ *—матрица* $A_{KN}$ *нормальна* ( $AA^* = A^*A$ ) *в том и только в том случае, если*

$$|\lambda_i(A_{KN})| = \sigma_i(A_{KN}), i = 0,1,2 \dots, N.$$

**Лемма 7.5.** *Пусть задано сингулярное разложение* $K \times N$ *—матрицы* $A_{KN} = U_{KK}\Sigma_{KN}V_{NN}^*$ *и* ${A_{KN}}^{\dagger}$ *— псевдообратная матрица. Тогда*

${A_{KN}}^{\dagger} = V_{NN}{A_{KN}}^{\dagger}{U_{KK}}^*.$

*Здесь*

$$ {A_{KN}}^{\dagger} = \operatorname{diag}\left(\sigma_1^{-1}, \sigma_2^{-2}, \cdots, \sigma_r^{-I}, \underbrace{0, \cdots, 0}_{\rho - r}\right), r = \operatorname{rank} A_{KN}, \rho = \min\{K, N\}, $$

т. е. ${A_{KN}}^{\dagger}$ — диагональная $K \times N$—матрица, у которой по диагонали стоят числа, обратные ненулевым сингулярным числам матрицы $A_{KN}$, а на остальных местах нули.

Еще одним замечательным свойством сингулярного разложения $K \times N$—матрицы $A_{KN} = U_{KK}\Sigma_{KN}V_{NN}^*$ является явное представление ядра и образа отображения $A_{KN}: \mathbb{R}^K \longrightarrow \mathbb{R}^N$, а именно, правые сингулярные векторы $v_{(k)}$ соответствующие нулевым сингулярным числам, порождают ядро $A_{KN}$, а левые сингулярные векторы $u_{(k)}$, соответствующие ненулевым сингулярным числам, порождают образ $R(A_{KN}) \subset \mathbb{R}^n$.

$$ U_{KK}\Sigma_{KN}V_{NN}^* q_N = f_N \tag{7.3} $$

$$ \Sigma_{KN} z = g \tag{7.4} $$

Пусть $K \times N$—матрица $A_{KN}$ имеет ранг $r \le \rho = \min\{K, N\}$. Тогда в сингулярном разложении $A_{KN} = U_{KK}\Sigma_{KN}V_{NN}^*$ матрица $\Sigma_{KN}$ имеет вид $\Sigma_{KN} = \operatorname{diag}(\sigma_\rho, \sigma_{\rho-1}, \cdots, \sigma_{\rho-r+1}, 0, \cdots, 0)$. При $r < K$ для совместности системы (7.4) необходимо, чтобы $g_{r+1} = g_{r+1} = \cdots = 0$. При $r < N$ переменные $z_{r+2}, z_{r+3}, \cdots, z_N$ в уравнение (7.4) не входят и, значит, в случае совместности системы компоненты решения $z_{r+1}, z_{r+2}, \cdots, z_N$ могут быть выбраны произвольно.

Выберем $z$ так, чтобы норма вектора. Невязки $\Sigma z - g = \begin{pmatrix} \sigma_\rho z_1 - g_1 \\ \sigma_{\rho-1} z_2 - g_2 \\ \vdots \\ \sigma_{\rho-r+1} z_r - g_r \\ -g_{r+1} \\ \vdots \\ -g_K \end{pmatrix}$, была минимальной (Годунов, 1988) [52].